\documentclass[11pt, reqno]{amsart}
\usepackage{amsmath, amssymb, amscd, amsthm, amsxtra}
\usepackage{graphicx}
\usepackage{color}
\usepackage{comment}
\usepackage{mathtools}
\usepackage[pdftex]{hyperref}

\usepackage[format=plain,labelfont=bf,up,width=.9\textwidth]{caption}

\theoremstyle{plain}
\newtheorem{thm}{Theorem}[section]
\newtheorem{cor}{Corollary}[thm]
\newtheorem{lemma}[thm]{Lemma}
\newtheorem{prop}[thm]{Proposition}
\newtheorem{remark}[thm]{Remark}
\theoremstyle{definition}

\theoremstyle{plain}
\newtheorem{thmx}{\bf Theorem}
\newtheorem{corx}{Corollary}[thmx]

\newcommand{\R}{\mathbb{R}}
\newcommand{\N}{\mathbb{N}}
\newcommand{\C}{\mathbb{C}}
\newcommand{\Z}{\mathbb{Z}}
\newcommand{\D}{\mathbb{D}}
\newcommand{\Q}{\mathbb{Q}}
\newcommand{\s}{\mathbb{S}}
\newcommand{\bigO}{\mathcal{O}}

\newcommand{\paratt}{\mathcal{P}_{\rm{par\mbox{-}att}}}

\newcommand{\Bparatt}{\mathcal{B}_{\rm{par\mbox{-}att}}}

\newcommand{\hog}{\mathcal{H}}
\newcommand{\Wloc}{W^c_{\rm loc}} 

\newcommand{\WlocO}{W^c_{\rm loc}(0)} 
\newcommand{\Chx}{\mathcal{C}_{x}^{h,\alpha}} 
\newcommand{\Cvx}{\mathcal{C}_{x}^{v,\alpha}} 
\newcommand{\ang}{{\rm Angle} }
\newcommand{\dbar}{\bar{\partial}}

\def\proof{\noindent {\bf Proof.\ }}

\def\qed{\hfill $\square$\\ }

\def \diag #1#2#3#4#5#6#7#8{\begin{CD} #1 @>#5>> #2\\ @V#6VV  @VV#7V\\ #3 @>#8>> #4 \end{CD}}
\def \diagup #1#2#3#4#5#6#7#8{\begin{CD} #1 @>#5>> #2\\ @A#6AA  @AA#7A\\ #3 @>#8>> #4 \end{CD}}

\begin{document}
\title[Hedgehogs in higher dimensions and their applications]{Hedgehogs in higher dimensions and their applications}
\author{Mikhail Lyubich}
\address{Stony Brook University, Stony Brook, United States}
\email{mlyubich@math.stonybrook.edu}

\author{Remus Radu}
\address{Stony Brook University, Stony Brook, United States}
\email{rradu@math.stonybrook.edu}

\author{Raluca Tanase}
\address{Stony Brook University, Stony Brook, United States}
\email{rtanase@math.stonybrook.edu}



\begin{abstract}
In this paper we study the dynamics of germs of holomorphic diffeomorphisms of $(\mathbb{C}^{n},0)$ with a fixed point at the origin with exactly one neutral eigenvalue. We prove that the map on any local center manifold of $0$ is quasiconformally conjugate to a holomorphic map and use this to transport results from one complex dimension to higher dimensions. 
\end{abstract}

\maketitle


\section{Introduction}\label{sec:Intro}

Let $f$ be a holomorphic germ of diffeomorphisms of $(\C^{2},0)$ with a fixed point at the origin with eigenvalues $\lambda$ and $\mu$, where $|\lambda|=1$ and $|\mu|<1$. Following the terminology from one-dimensional dynamics, the fixed point is called {\it semi-neutral} or {\it semi-indifferent}.

The crude analysis of the local dynamics of the semi-indifferent fixed point exhibits the existence of an analytic strong stable manifold $W^{ss}(0)$ corresponding to the dissipative eigenvalue $\mu$ and a not necessarily unique local center manifold $\WlocO$ corresponding to the neutral eigenvalue $\lambda$. 
The center manifolds can be made $C^r$-smooth for any $r\geq 1$ by possibly restricting to smaller neighborhoods of the origin \cite{HPS}, however they are generally not analytic, or even $C^{\infty}$-smooth \cite{vS}.
In this paper, we show how to modify the complex structure on the center manifold $\Wloc(0)$ so that the restriction of the map $f$ to the center manifold becomes analytic. 

\begin{thmx}\label{thm:QC} Let $f$ be a holomorphic germ of diffeomorphisms of $(\C^2,0)$ with a semi-neutral fixed point at the origin with eigenvalues $\lambda$ and $\mu$, where $|\lambda|=1$ and $|\mu|<1$. Consider $\WlocO$ a $C^1$-smooth local center manifold of the fixed point $0$. There exist neighborhoods $W, W'$ of the origin inside $\WlocO$ such that $f:W\rightarrow W'$ is quasiconformally conjugate to a holomorphic diffeomorphism $h:(\Omega,0)\rightarrow (\Omega',0)$, $h(z) = \lambda z + \bigO(z^2)$, where $\Omega, \Omega'\subset\C$.

The conjugacy map is holomorphic on the interior of $\Lambda$ rel $\WlocO$, where $\Lambda$ is the set of points that stay in $W$ under all backward iterations of $f$.
\end{thmx}

Theorem \ref{thm:QC} generalizes to the case of holomorphic germs of diffeomorphisms of $(\C^n,0)$, for $n>2$, which have a fixed point at the origin with exactly one eigenvalue on the unit circle. The details are given in Section \ref{sec:Cn}.

In dimension one, linearization properties and dynamics of  holomorphic univalent germs of $(\C,0)$ with an indifferent fixed point at the origin have been extensively studied (\cite{Y1}, \cite{Y2}, \cite{PM1}, \cite{PM2}, \cite{PM3}, and many more). Theorem \ref{thm:QC} has important consequences and enables to us to transport results from one complex variable to $\C^2$. In Section \ref{sec:PM} we examine the results of P\'erez-Marco about the hedgehog dynamics, and in Section  \ref{sec:cons} we show how to extend them to $\C^2$ using Theorem \ref{thm:QC}.

Suppose that the neutral eigenvalue $\lambda$ of the semi-neutral fixed point of the germ $f$ is $\lambda=e^{2\pi i \alpha}$. If the origin is an isolated fixed point of $f$ and $\alpha\in\Q$, 
then the fixed point is called {\it semi-parabolic}. In the case when $\alpha\notin\Q$, the fixed point is called {\it irrational semi-indifferent}. We can further classify irrational semi-indifferent fixed points as semi-Siegel or semi-Cremer, as follows:  if there exists an injective holomorphic map $\varphi:\D\rightarrow \C^2$ such that $f(\varphi(\xi))=\varphi(\lambda \xi)$, for $\xi\in\D$, then the fixed point is called {\it semi-Siegel}, otherwise it is called {\it semi-Cremer}. 
Theorem \ref{thm:D} below motivates the following equivalent definition: if $f$ is analytically conjugate to $(x,y)\mapsto(\lambda x, \mu(x)y)$, where $\mu(x)=\mu+\bigO(x^2)$ is a holomorphic function, then the fixed point is {\it semi-Siegel}; otherwise, the fixed point is {\it semi-Cremer}. In particular, when $\lambda$ satisfies the Brjuno condition \cite{Brj} and $|\mu|<1$, the map $f$ is linearizable ({\it i.e.} conjugate by a holomorphic change of variables to its linear part), so the fixed point is semi-Siegel.

In \cite{FLRT} we have shown the existence of non-trivial compact invariant sets for germs $f$ with semi-indifferent fixed points, using topological tools. 
If $f$ is a germ of holomorphic diffeomorphisms of $(\C^{2},0)$
with a semi-indifferent fixed point at the origin, then there exists a domain $B$ containing $0$ such that $f$ is partially hyperbolic on a neighborhood of $\overline{B}$. The concept of partial hyperbolicity is explained in the introductory part of Section \ref{sec:structure}.

\begin{thm}[\cite{FLRT}]\label{thm:FLRT} 
Let $f$ be a germ of holomorphic diffeomorphisms of $(\C^{2},0)$
with a semi-indifferent fixed point at $0$ with eigenvalues $\lambda$ and $\mu$, where $|\lambda|=1$ and $|\mu|<1$.   
Consider an open ball $B\subset \C^{2}$ containing $0$ such that $f$ is partially hyperbolic on a neighborhood $B'$ of $\overline{B}$.
There exists a set $\mathcal{H}\subset \overline{B}$ such that:
\begin{itemize}
\item[a)] $\mathcal{H}\Subset W^{c}_{\rm loc}(0)$, where $W^{c}_{\rm loc}(0)$ is any local center manifold of the fixed point $0$, constructed relative to $B'$.
\item[b)] $\mathcal{H}$ is compact, connected, completely invariant, and full.
\item[c)] $0\in \mathcal{H}$, $\mathcal{H}\cap\partial B\neq \emptyset$. 
\item[d)] Every point $x\in\mathcal{H}$ has a well defined local strong stable manifold $W^{ss}_{\rm loc}(x)$, consisting of points from $B$ that converge asymptotically exponentially fast to $x$, at a rate $\asymp$ $\mu^{n}$. The strong stable set of $\mathcal{H}$ is laminated by vertical-like holomorphic disks.
\end{itemize}
\end{thm}

In this paper, the compact set $\mathcal{H}$ will be called a {\it hedgehog}. We distinguish between a {\it parabolic hedgehog}, a {\it Siegel hedgehog}, or a {\it Cremer hedgehog} (also called {\it non-linearizable hedgehog}), depending whether the fixed point is semi-parabolic, semi-Siegel, or semi-Cremer.  
In this paper we will explore the dynamical properties of hedgehogs. 
Let us start by noticing that in the irrationally semi-indifferent case the hedgehog is locally unique. There exists a neighborhood $B$ of the origin in $\C^{2}$ such that the hedgehog $\mathcal{H}$  associated to $f$ and $B$ is unique and equal to the connected component containing $0$ of the set $\{z\in \overline{B} : f^n(z) \in \overline{B}\ \forall n\in \Z\}$. This is a direct consequence of Theorem \ref{thm:QC} and Theorem \ref{thm:PM2} of P\'erez-Marco.

\begin{figure}[htb]
\begin{center}
\includegraphics[scale=.59]{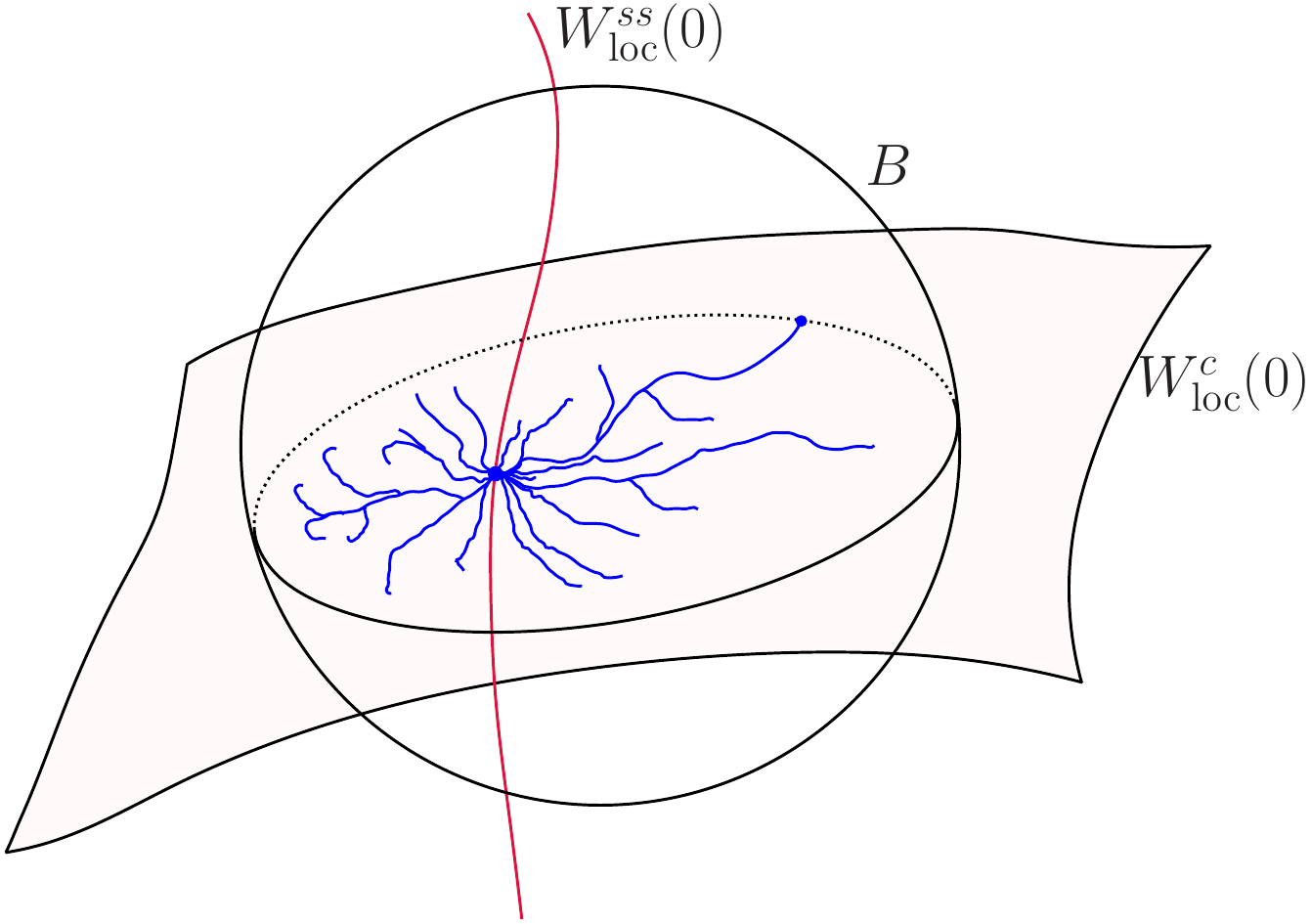}
\end{center}
\caption{A Cremer hedgehog inside a center manifold.}
\label{fig:hedgehog}
\end{figure}

The next two theorems and the subsequent corollaries deal with Cremer hedgehogs (see Figure \ref{fig:hedgehog}).

\begin{thmx}\label{thm:B} Let $f$ be a germ of dissipative holomorphic diffeomorphisms of $(\C^{2},0)$ with a semi-Cremer fixed point at $0$ with an eigenvalue $\lambda=e^{2\pi i \alpha}$. Let $\hog$ be a hedgehog for $f$. Suppose $(p_{n}/q_{n})_{n\geq 1}$ are the convergents of the continued fraction of $\alpha$. There exists a subsequence $(n_k)_{k\geq 1}$ such that the iterates $(f^{q_{n_k}})_{k\geq 1}$ converge uniformly on $\hog$ to the identity.  
\end{thmx}

\begin{corx}\label{corB1}
The dynamics on the hedgehog $\hog$ is recurrent. The hedgehog does not contain other periodic points except $0$. 
\end{corx}

Denote by $\omega(x)$ and $\alpha(x)$ the $\omega$-limit, respectively $\alpha$-limit set of $x$. 

\begin{thmx}\label{thm:C} There exists a neighborhood $B\subset \C^2$ of the origin with the following properties. Let $\mathcal{H}$ be a Cremer hedgehog for $f$ and $B$.  

\begin{itemize}
\item[a)] Let $x\in B-W_{\rm loc}^{ss}(\mathcal{H})$. If the forward iterates $f^{n}(x)\in B$ for all $n\geq 0$, then $\omega(x)\cap \mathcal{H}=\emptyset$.
\item[b)] Let $x\in B-\mathcal{H}$.  If the backward iterates $f^{-n}(x)\in B$ for all $n\geq 0$, then $\alpha(x)\cap \mathcal{H}=\emptyset$.
\end{itemize}
\end{thmx}

Theorem \ref{thm:C} immediately implies the following corollaries. 
\begin{corx}\label{corC1}
If $x\notin W^{ss}(0)$ then the orbit of $x$ does not converge to $0$. 
\end{corx}

Let $\mathcal{H}$ be a hedgehog for a germ $f$ with a semi-Cremer fixed point and denote by
 $W^{s}_{\rm loc}(\mathcal{H})$ the local stable set of $\mathcal{H}$, consisting of points that converge to the hedgehog. Let $W^{ss}_{\rm loc}(\mathcal{H})$ be the local strong stable stable of $\mathcal{H}$, consisting of points which converge asymptotically exponentially fast to the hedgehog. 
\begin{corx}\label{corC2}
$ W^{ss}_{\rm loc}(\mathcal{H})=W^{s}_{\rm loc}(\mathcal{H})$.
 \end{corx}

Clearly, the set $\mathcal{H}$ has no interior in $\C^2$ since it lives in a center manifold. Let ${\rm int}^{c}(\mathcal{H})$ denote the interior of $\mathcal{H}$ relative to a center manifold. 

\begin{thmx}\label{thm:D} Let $f$ be a holomorphic germ of diffeomorphisms of $(\C^2,0)$ with an isolated semi-neutral fixed point at the origin. Let $\mathcal{H}$ be a hedgehog for $f$. Then $0\in {\rm int}^{c}(\mathcal{H})$ if and only if $f$ is analytically conjugate to a linear cocycle $\tilde{f}$ given by 
\[
\tilde{f} (x,y) = (\lambda x, \mu(x)y),
\]
where $\mu(x) = \mu+\bigO(x)$ is a holomorphic function. 
\end{thmx}

\begin{corx}\label{cor:lin} Let $f$ be a dissipative polynomial diffeomorphism of $\C^2$ with a semi-neutral fixed point at the origin. Then $0\in {\rm int}^{c}(\mathcal{H})$ if and only if $f$ is linearizable.
\end{corx}

For a polynomial automorphism $f$ of $\C^2$ we define the sets $K^{\pm}$ of points that do not escape to $\infty$ under forward/backward iterations. Denote by $J^{\pm}$ the topological boundaries of $K^{\pm}$ in $\C^2$. The set $J=J^-\cap J^+$ is called the Julia set. 
Let $J^{*}$ be the closure of the saddle periodic points of $f$. The set $J^{*}$ is a subset of $J$, but it is an open question whether they are always equal.
We obtain the following result.

\begin{thmx}\label{thm:wandering} Let $f$ be a dissipative polynomial diffeomorphism of $\C^2$ with an irrationally semi-indifferent fixed point at $0$. Suppose $f$ is not linearizable in a neighborhood of the origin. Let $\mathcal{H}$ be a hedgehog for $f$. Then $\mathcal{H}\subset J^{*}$ and there are no wandering domains converging to $\mathcal{H}$.
\end{thmx}

Consider now a holomorphic germ $f$ of $(\C^2,0)$,  with a semi-parabolic fixed point at $0$ with a neutral eigenvalue 
$\lambda=e^{2\pi i p/q}$. After a holomorphic change of coordinates, we can assume that $f$ is written in the normal form 
$f(x,y)=(x_1,y_1)$, where
\begin{equation*}\label{eq:NF}
\left\{\begin{array}{l}
    x_{1}= \lambda (x+x^{\nu q+1} + Cx^{2\nu q+1}+ a_{2\nu q+2}(y)x^{2\nu q +2}+\ldots )\\
    y_{1}= \mu y + xh(x,y)
\end{array}\right.
\end{equation*}
We call $\nu$ the {\it semi-parabolic multiplicity} of the semi-parabolic fixed point.

The existence of holomorphic 1-D repelling petals for $f$  and of Fatou coordinates on the repelling petals, was established by Ueda  \cite{U2}.  We can give a new proof of this result using Theorem \ref{thm:QC} and the Leau-Fatou theory of parabolic holomorphic germs of $(\C,0)$. Let $B\subset\C^2$ be a small enough ball around the origin and define the set 
\begin{equation}\label{eq:Sigma}
\Sigma_B=\{x\in B-\{0\} : f^{-n}(x)\in B\ \forall n\in\N,\ f^{-n}(x)\rightarrow 0 \mbox{ as } n\rightarrow\infty\}.
\end{equation}

\begin{thmx}\label{thm:F} Let $f$ be a  germ of dissipative holomorphic diffeomorphisms of $(\C^2,0)$ with a  semi-parabolic fixed point at $0$ with an eigenvalue $\lambda=e^{2\pi i p/q}$. Suppose the semi-parabolic multiplicity of $0$ is $\nu$. 
The set $\Sigma_B$ from \eqref{eq:Sigma} consists of $\nu$ cycles of $q$ repelling petals. Each repelling petal 
is the image of an injective holomorphic map $\varphi(x)=\left(x,k(x)\right)$ from a left half plane of $\C$ into $\C^2$, which satisfies $\varphi (x+1)=f^{q}(\varphi(x))$. The inverse of $\varphi$, denoted by $\varphi^{o}:\Sigma_B\rightarrow \C$ is called an outgoing Fatou coordinate; it satisfies the Abel equation $\varphi^{o}(f^q)=\varphi^{o}+1$.
\end{thmx}

\noindent\textit{Acknowledgement.} We thank Romain Dujardin for pointing out that in Theorem \ref{thm:wandering} the hedgehog $\mathcal{H}$ is in $J^{*}$, not just in  $J$. The first author was supported by NSF grants DMS-1301602 and DMS-1600519.

\section{Holomorphic germs of $(\C,0)$}\label{sec:PM}
In this section we make a brief survey of the results of P\'erez-Marco on holomorphic germs of $(\C,0)$ with a neutral fixed point at the origin. The following theorem gives the local structure of such germs.

\begin{thm}[\cite{PM1}]\label{thm:PM1} Let $f(z)=\lambda z +\bigO(z^2)$, $|\lambda|=1$ be a local holomorphic diffeomorphism, and $U$ a Jordan neighborhood of the indifferent fixed point $0$. Assume that $f$ and $f^{-1}$ are defined and univalent on a neighborhood of the closure of $U$. Then there exists a completely invariant set $K\subset \overline{U}$, compact connected and full, such that $0\in K$ and $K \cap \partial U \neq \emptyset$. 

Moreover, if $f^{ n}\neq id$ for all $n\in \N$ then $f$ is linearizable if and only if $0\in {\rm int}(K)$. 
\end{thm}

P\'erez-Marco calls the compact set $K$ a Siegel compactum for $f$. He calls $K$ a hedgehog when the fixed point is irrationally indifferent, and the set $K$ is not contained in the closure of a linearization domain. For simplicity, we refer to the sets $K$ from Theorem \ref{thm:PM1} as {\it hedgehogs} and we distinguish between the various types of hedgehogs as in the introduction.

Let $K$ be a hedgehog for $f$, as in Theorem \ref{thm:PM1}, and $\lambda=e^{2\pi i \alpha}$. P\'erez-Marco associates to each set $K$ an analytic circle diffeomorphism with rotation number $\alpha$ as follows. We first uniformize $\C\setminus K$ using the Riemann map $\psi:\C\setminus K\rightarrow \C\setminus \overline{\D}$. Let  $g=\psi\circ f\circ \psi^{-1}$. The mapping $g$ is defined and holomorphic in an open annulus $\{z\in \C: 1<|z|<r\}$. We can extend $g$ to the annulus $\{z\in \C: 1/r<|z|<r\}$ by the Schwarz reflexion principle. The restriction $g|_{\s^1}$ to the unit circle will be a real-analytic diffeomorphism with rotation number $\alpha$. P\'erez-Marco uses this construction to transport results of analytic circle diffeomorphisms to results about the dynamics of holomorphic maps around the indifferent fixed points.  This analogy is used both to construct the hedgehog from Theorem \ref{thm:PM1} and  to study the dynamics on the hedgehog. 

We say that a domain $U$ is {\it admissible} if it is a $C^1$-Jordan domain such that $f$ and $f^{-1}$ are univalent on a neighborhood of the closure of $U$. 

 \begin{thm}[\cite{PM1}, \cite{PM3}]\label{thm:PM2} Let $U$ be an admissible neighborhood for a germ $f(z)=\lambda z +\bigO(z^2)$, $\lambda=e^{2\pi i \alpha}$ and $\alpha\notin\Q$.
 Let $L$ be a connected compact set, invariant under $f$ or $f^{-1}$, such that $0\in L\subset \overline{U}$ and $L\cap \partial U\neq\emptyset$. Then $L=K$. 
 Therefore the hedgehog $K$ of $f$ associated to the neighborhood $U$ is unique and it is equal to the connected component
containing $0$ of the set $\{z\in \overline{U} : f^n(z) \in \overline{U} \ \mbox{for all}\ n\in \Z\}$. 
 \end{thm}

For the rest of the section, suppose that $f$ is a non-linearizable germ with $f(z)=\lambda z +\bigO(z^2)$, $\lambda=e^{2\pi i \alpha}$ and $\alpha\notin\Q$. Consider $U$ an admissible neighborhood and $K$ a hedgehog associated to $f$ and $U$. Let $(p_{n}/q_{n})_{n\geq 1}$ be the convergents of the continued fraction of $\alpha$. P\'erez-Marco shows that even if $f$ is non-linearizable, the dynamics of $f$ has a lot of features in common to the irrational rotation.

\begin{thm}[\cite{PM2},\cite{PM3}]\label{thm:PM3}
There exists a subsequence $A\subset\N$ such that the iterates $(f^{q_{n}})_{n\in A}$ converge uniformly on $K$ to the identity. 
\end{thm}

\begin{cor}\label{cor:PM3} All points on the hedgehog $K$ are recurrent. The dynamics on the hedgehog has no periodic point except the fixed point at $0$. 
\end{cor}

\begin{thm}[\cite{PM2},\cite{PM3}]\label{thm:PM4} Let $x\in U$ be a point which does not belong to $K$. If the $\omega$-limit (or $\alpha$-limit) set of $x$ intersects $K$, then it cannot be contained in $U$. Consequently, $x$ cannot converge to a point of $K$ under iterations of $f$. 
\end{thm}

Regarding the topology of the hedgehog,  P\'erez-Marco \cite{PM3} shows that non-linearizable hedgehogs have empty interior and are not locally connected at any point different from $0$. 

The following two theorems discuss the structure of non-linearizable germs.  Using Theorem \ref{thm:QC}, one can immediately formulate a higher dimensional analogue of the following result:

\begin{thm}[\cite{PM4}, \cite{PM5}]\label{thm:PM4-1} If 
\begin{equation}\label{eq:loglog}
	\sum\limits_{n=1}^{\infty} \frac{\log\log q_{n+1}}{q_n} < \infty
\end{equation}
then all non-linearizable germs $f(z) = \lambda z + \bigO(z^2)$ have a sequence of periodic orbits $(O_k)_{k\geq 0}$ which tend to $0$, of periods $q_{n_k}$ and rotation numbers $p_{n_k}/q_{n_k}$ such that 
\begin{equation*}
	\sum\limits_{n=1}^{\infty} \frac{\log q_{n_{k}+1}}{q_{n_{k}}} = \infty.
\end{equation*}
\end{thm}

Condition \eqref{eq:loglog} is sharp:
by \cite{PM5}, if the sum in \eqref{eq:loglog} diverges, then there exists a holomorphic germ $f(z) = \lambda z + \bigO(z^2)$ defined and univalent in $\D$ with no other periodic orbits in $\D$, except $0$. In fact, every orbit $(f^n(z))_{n\geq0}$ remaining in $\D$ accumulates $0$, {\it i.e.} $0\in \omega(z)$.

\section{A complex structure on the center manifold}\label{sec:structure}

Consider a holomorphic germ $f$ of diffeomorphisms of $(\C^2,0)$ with a semi-indifferent fixed point at $0$. 
Denote by $\lambda$ and $\mu$ the eigenvalues of the derivative $df_0$. Throughout the paper, we make the convention that $|\lambda|=1$ (neutral eigenvalue) and $|\mu|<1$ (dissipative eigenvalue), and we denote by $E_0^c$ and $E_0^s$ the eigenspaces corresponding to the eigenvalues $\lambda$ and respectively $\mu$. The map $f$ is partially hyperbolic on a neighborhood $B'\subset \C^2$ of the origin. 
Partial hyperbolicity was introduced in the '70s by Brin and Pesin \cite{BP}  and Hirsch, Pugh, and Shub \cite{HPS} as a natural extension of the concept of hyperbolicity. 
We define partial hyperbolicity below; for a thorough introduction and an overview of the field we refer the reader to \cite{HP}, \cite{CP}, and \cite{HHU}.

Let $E^1$ and $E^2$ be two continuous distributions (not necessarily invariant by $df$) of the complex tangent bundle $TB'$ such that $T_x B'=E^1_x\oplus E^2_x$ for any $x\in B'$. For $x=0$, we take $E_0^1=E_0^c$ and $E_0^2=E_0^s$. The horizontal cone $\Chx$ is defined as the set of vectors in the tangent space at $x$ that make an angle less than or equal to $\alpha$ with $E^{1}_x$, for some $\alpha>0$,
\[
\mathcal{C}^{h,\alpha}_{x}=\{ v\in T_x B', \angle(v,E_{x}^{1})\leq \alpha \},
\]
where the angle of a vector $v$ and a subspace $E$ is simply the angle between $v$ and its projection $\mbox{pr}_E v$ on the subspace $E$.
The vertical cone $\Cvx$ is defined in the same way, with respect to $E^{2}_x$. We will suppress the angles $\alpha$ from the notation of the cones, whenever there is no danger of confusion. 

The map $f$ is partially hyperbolic on $B'$ if there exist two real numbers $\mu_1$ and $\lambda_1$ such that $0<|\mu|<\mu_1<\lambda_1<1$ and a family of invariant cone fields $\mathcal{C}^{h/v}$
\begin{equation}\label{eq:inv}
df_x(\mathcal{C}_x^{h})\subset \mbox{ Int } \mathcal{C}_{f(x)}^{h}\cup\{0\}, \ \ \ df^{-1}_{f(x)}(\mathcal{C}_{f(x)}^{v})\subset \mbox{ Int } \mathcal{C}_{x}^{v}\cup\{0\},
\end{equation}
such that for some Riemannian metric we have strong contraction in the vertical cones, whereas in the horizontal cones we may have contraction or expansion, but with smaller rates:
\begin{equation}\label{eq:Cn}
\lambda_1\,\|v\|\leq \|df_x(v)\| \leq \lambda_1^{-1}\,\|v\|,\ \ \ \ \ \mbox{for}\ v\in \mathcal{C}_x^{h}
 \end{equation}
\begin{equation*}\label{eq:Cs}
   \|Df_x(v)\|\leq \mu_{1}\, \|v\|,\ \ \ \ \ \mbox{for}\ v\in \mathcal{C}_{x}^{v}.
\end{equation*}

Let $B$ be a domain in $\C^{2}$ containing the origin such that $\overline{B}\subset B'$ and $f(\overline{B})\subset B'$. 

\begin{remark}\label{rem:rho} Since $B$ is compactly contained in $B'$, 
the angle between $df_x(\mathcal{C}_x^{h})$ and $\partial \mathcal{C}_{f(x)}^{h}$ is uniformly bounded independently of $x$. This implies that there exists 
$0<\rho<1$  such that for every $x\in B$, the angle opening of the cone $df_x(\mathcal{C}_x^{h})$ is $\rho\alpha$, a fraction of the angle opening of the cone $\mathcal{C}_{f(x)}^{h}$.
\end{remark}

 The semi-indifferent fixed point has a local center manifold $\Wloc:=\WlocO$ of class $C^1$, tangent at $0$ to the eigenspace $E_0^c$ corresponding to the neutral eigenvalue $\lambda$. Throughout the paper, $\Wloc$ will denote the local center manifold of $0$. 
The local center manifold is the graph of a $C^1$ function $\varphi_f: E_0^c \cap B' \rightarrow E_0^s$ and  has the following properties: 
\begin{itemize}
\item[a)] {\it Local invariance} $f(\Wloc)\cap B'\subset \Wloc$.
\item[b)] {\it Weak uniqueness:} 
If $f^{-n}(x)\in B'$ for all $n\in\N$, then $x\in \Wloc$. 
\item[c)] {\it Shadowing:} Given any point $x$ such that $f^{n}(x)\rightarrow 0$ as $n\rightarrow \infty$, there exists a positive constant $k$ and a point $y\in \Wloc$ such that $\|f^{n}(x)-f^{n}(y)\|<k\mu_1^n$ as $n\rightarrow \infty$. In other words, every orbit which converges to the origin can be described as an exponentially small perturbation of some orbit on the center manifold.  
\end{itemize} 

The center manifold is generally not unique. However, the formal Taylor expansion at the origin is the same for all center manifolds. The center manifold is unique in some  cases, for instance when $f$ is complex linearizable at the origin. For uniqueness, existence, and regularity properties of center manifolds, we refer the reader to \cite{HPS}, \cite{Sij}, \cite{S}, and \cite{V}. 

It is also worth mentioning the following {\it reduction principle} for center manifolds: the map $f$ is locally topologically semi-conjugate to a function on the center manifold given by $u\mapsto \lambda u + f_1(u,\varphi_{f}(u))$. In this article, we will not make use of the reduction principle, as it only gives a topological semi-conjugacy to a model map which is as regular as the center manifold, hence not analytic.

The assumption of partial hyperbolicity implies that any point $x$ with $f^{n}(x)\in B$ for all $n\geq 0$ has a well defined local strong stable manifold $W^{ss}_{\rm loc}(x)$, given by 
\[
W^{ss}_{\rm loc}(x) = \{y\in B : {\rm dist}(f^n(x),f^n(y))/\mu_1^n\rightarrow 0\ \mbox{as}\ n\rightarrow \infty\}.
\]
Moreover $W^{ss}_{\rm loc}(x)$  intersects $\Wloc$ transversely. Let $y\in W^{ss}_{\rm loc}(x)\cap \Wloc$. Then the orbit of $y$ shadows the orbit of $x$. We can therefore formulate a more general shadowing property:

\begin{prop}\label{prop:shadow}
Let $x\in B$ such that $f^n(x)\in B$ for every $n\geq 1$. 
There exists $k>0$ and $y\in \Wloc$ such that $\|f^{n}(x)-f^{n}(y)\|<k\mu_1^n$ as $n\rightarrow \infty$. 
\end{prop}

An old question posed by Dulac and Fatou is whether there exists orbits converging to an irrationally indifferent fixed point of a holomorphic map. Using the dynamics on the hedgehogs, renormalization theory and Yoccoz estimates for analytic circle diffeomorphisms, P\'erez-Marco showed that the answer to this question is negative, see Theorem \ref{thm:PM4}.  It is a natural question to ask whether there exist orbits converging to a semi-Cremer fixed point of a holomorphic germ of $(\C^2,0)$. 
Also, P\'erez-Marco has constructed examples of hedgehogs in which the origin is accumulated by periodic orbits of high periods, see Theorem \ref{thm:PM4-1}.

In the two-dimensional setting, the problem of the existence of periodic orbits of $f$ accumulating the origin or the existence of orbits converging to zero can naturally be reduced to posing the same question for the restriction of the map $f$ to the local center manifold(s). 
Let us explain this reduction further.

If the semi-indifferent fixed point is accumulated by periodic orbits of high period, then these periodic points necessarily live in the intersection of all center manifolds $\WlocO$, by the weak uniqueness property of local center manifolds. 
Since we work with dissipative maps, there will always be points that converge to $0$ under forward iterations, corresponding to the strong stable manifold of $0$. If there exists some other point $x$, whose forward orbit converges to $0$, then the orbit of $x$ must be shadowed by the orbit of a point $y$ that converges to $0$ on the center manifold. Therefore, in order to answer the questions about the dynamics of the two-dimensional germ around $0$, we should first study the dynamics of the function restricted to the center manifolds.

The main obstruction for extending the results of P\'erez-Marco directly to our setting is the fact that the center manifolds are not analytic. It is well-known (see e.g. \cite{GH}, \cite{vS}) that there exist $C^{\infty}$-smooth germs which do not have any $C^{\infty}$-smooth center manifolds. For every finite $k$ one can find a neighborhood $B_k$ of the origin for which there exists a $C^k$-smooth center manifold relative to $B_k$, however the sets $B_k$ shrink to $0$ as $k\rightarrow \infty$.  

A first useful observation in this context is that some analytic structure still exists in some parts of the center manifold. 
Let 
\begin{equation}\label{eq:Lambda}
\Lambda = \{z\in B : f^{-n}(z)\in B,\ \mbox{for all}\ n\geq 0\}
\end{equation}
be the set of points that never leave $B$ under backward iterations. 
As a consequence of Theorem \ref{thm:FLRT}, we know that the set $\Lambda$ is not trivial, {\it i.e.} $\Lambda\neq\{0\}$.  By the weak uniqueness property of center manifolds, $\Lambda$ is a subset of $\Wloc$. Also, $f^{-1}(\Lambda)\subset \Lambda$, by definition.

\begin{prop}\label{horiz} 
Let $\Wloc$ be any local center manifold relative to $B'$. The tangent space $T_{x}\Wloc$ at any point $x\in \Lambda$ is a complex line $E^c_x$ of $T_x\C^2$. The line field over $\Lambda$ is $df$-invariant, in the sense that $df_x(E^c_x)=E^c_{f(x)}$ for every point $x\in\Lambda$ with $f(x)\in\Lambda$.
\end{prop}
\proof
Let $x\in  \Lambda$. All iterates $f^{-n}(x)$, $n\geq 0$  remain in the domain $B$, where we have an invariant family of horizontal cones $\mathcal{C}^h$ preserved by $df$. The derivative acts on tangent vectors as a vertical contraction by a factor $\mu_1$, close to $\mu$. Let 
\[
E^c_x= \bigcap_{n\geq 0} df^n_{f^{-n}(x)} \mathcal{C}^h_{f^{-n}(x)}.
\] 
This is a decreasing intersection of nontrivial compact subsets in the projective space, hence $E^c_x$ is a non-trivial complex subspace of $T_{x}\C^2$. Thus
$E^c_x$ is a complex line included in $\mathcal{C}^h_x$. The invariance of the line bundle $(E^c_x)_{x\in\Lambda}$ under $df$ follows from the definition. 

Any local center manifold $\Wloc$ defined relative to $B'$  contains the set $\Lambda$ and the tangent space $T_{x}\Wloc$ at any point $x\in \Lambda$ is equal to $E^c_x$, thus it is a complex line in the tangent bundle $T\C^2$.
\qed

Proposition \ref{horiz}  means in particular that for every point $x\in \Lambda$, the tangent space $T_x(\Wloc)$ is $J$-invariant, where $J$ is the standard almost complex structure obtained from the usual identification of $\R^{4}$ with $\C^2$. Recall that an almost complex structure on a smooth even dimensional manifold $M$ is a complex structure on its tangent bundle $TM$, or equivalently a smooth $\R$-linear bundle map $J:TM\rightarrow TM$ with $J\circ J=-Id$.

The center manifold $\Wloc$ is a real 2-dimensional submanifold of $\C^2$. The standard Hermitian metric of the complex manifold $\C^2$ defines a Riemannian metric on the underlying smooth manifold $\R^4$, which restricts to a Riemannian metric on the center manifold $\Wloc$. Recall that in $\C^n$, the standard Hermitian inner product decomposes into its real and imaginary parts:
$
\langle u,v\rangle_H = \langle u, v\rangle-iw(u,v),
$
where $\langle u, v\rangle$ is the Euclidean scalar product and $w(u,v)$ is the standard symplectic form of $\R^{2n}$. From now on, whenever we refer to the  Riemannian metric we understand the metric defined by the Euclidean scalar product.

Every Riemannian metric on an oriented 2-dimensional manifold induces an almost complex structure given by the rotation by $90^{\circ}$, {\it i.e.} by defining  
\[
J'_x:T_x\Wloc\rightarrow T_x\Wloc,\ \ \mbox{as}\ \ J'_x(v)=v^{\perp},
\] 
where $v^{\perp}$ is the unique vector orthogonal to $v$, of  norm equal to $\|v\|$,
 such that the choice is orientation preserving. Every almost complex structure on a 2-dimensional manifold is integrable, that is, it arises from an underlying complex structure. Namely, there exists a $(J',i)$-holomorphic parametrization function $\phi : \Delta \rightarrow \Wloc$ where $\Delta$ is an open subset in $\C$ and $i$ is the standard complex structure in $\C$ given by multiplication by the complex number $i$. By $(J',i)$-holomorphic map, we understand a $C^1$-smooth map with the property that its derivative $d\phi_z:T_z\Delta\rightarrow T_{\phi(z)}\Wloc$ is complex linear, that is $d\phi_z \circ i_z = J'_{\phi(z)}\circ d\phi_z$. A good introduction on almost complex structures and $J$-holomorphic curves can be found in \cite{Voi} and \cite{MS}. Note that the parametrizing map $\phi$ that we have constructed is only $(J',i)$-holomorphic, but not necessarily $(J,i)$-holomorphic, as $\Wloc$ is not in general an embedded complex submanifold of $\C^2$. Note also that the almost complex structure induced by $J'$ on $\Wloc$ agrees with the standard almost complex structure $J$ from $\C^2$ on the set $\Lambda$, so $J_x=J'_x$ for all $x\in \Lambda$. 

Let $W := B\cap \Wloc$ and $U:=\phi^{-1}(W)\subset \Delta$. The set $W'=f(W)$ belongs to $\Wloc$, by the local invariance of the center manifold.  The map $f:W\rightarrow W'$ is an orientation-preserving $C^1$ diffeomorphism. 

Let $g=\phi\circ f\circ \phi^{-1}:U\rightarrow U'=g(U)$ be the  orientation-preserving $C^1$-diffeomorphism induced by $f$ on $U$. 

\[
\diagup{W}{W'}{U}{U'}{f}{\phi}{\phi}{g}
\]

\medskip
Denote by $X:=\phi^{-1}(\Lambda)$, or equivalently
\begin{equation}\label{eq:X}
X= \{z\in U : g^{-n}(z)\in U,\ \mbox{for all}\ n\geq 0\},
\end{equation}
the set of points that stay in $U$ under all backward iterations by $g$. 

 The map $f$ is holomorphic on $\C^2$, so it is $(J,J)$-holomorphic on $\Lambda$, which means that  $\bar{\partial}_J f =0$ for $\xi \in \Lambda$, where
\begin{equation}\label{eq:dbar1}
\bar{\partial}_J f := \frac{1}{2}(df_{\xi}+J_{f(\xi)}\circ df_{\xi} \circ J_{\xi}).
\end{equation} 
The conjugacy function $\phi$ is $(J,i)$-holomorphic on $\Lambda$. In the holomorphic coordinates provided by $\phi$, this means that $g$ is $(i,i)$-holomorphic on $X$, or equivalently $\bar{\partial}_i g=0$, where
\begin{equation}\label{eq:dbar}
\bar{\partial}_i g := \frac{1}{2}(dg_z+i_{g(z)}\circ dg_z \circ i_z),
\end{equation} 
and $z=\phi^{-1}(\xi)$. 

It is easy to check that with the standard identifications $z=x+iy$, $g(x,y)=g_1(x,y)+ig_2(x,y)$, and 
\begin{equation}\label{eq:iz}
i_z\left(\partial_x\right)=\partial_y, \ \ \  i_z\left(\partial_y\right)=-\partial_x,
\end{equation}
the relation $\bar{\partial}_i g=0$ is equivalent to the familiar Cauchy-Riemann equations 
\[
\partial_x g_1-\partial_y g_2=0\ \ \mbox{and}\ \ \partial_x g_2+\partial_y g_1=0.
\]

As usual, consider the linear partial differential operators of first order
\begin{equation}\label{eq:dd}
\partial = \partial_z = \frac{1}{2}\left(\partial_x-i \partial_y \right)\ \ \ \mbox{and}\ \ \ \bar{\partial} = \partial_{\bar{z}}= \frac{1}{2}\left(\partial_x+i\partial_y  \right).
\end{equation}
We have just shown the following proposition:
\begin{prop}\label{prop:dbar-X}  
$\bar{\partial} g = 0$ on the set $X$,  defined in Equation \eqref{eq:X}.
\end{prop}

Let ${\rm int}^c(\Lambda)$ denote the interior of $\Lambda$ rel $\Wloc$. Propositions \ref{horiz} and \ref{prop:dbar-X} immediately imply the following corollary:

\begin{cor}\label{cor:submanC} The set ${\rm int}^c(\Lambda)$ is a complex submanifold of $\C^2$. The conjugacy map $\phi:{\rm int}(X)\subset\C\rightarrow {\rm int}^c(\Lambda)\subset\C^2$ is holomorphic, and the function $g$ is holomorphic on ${\rm int}(X)$.
\end{cor}

When $x\in\Wloc-\Lambda$, the tangent space $T_x\Wloc$ is only a real 2-dimensional subspace of $T_x\C^2$, included in the horizontal cone $\mathcal{C}^h_x$. We want to measure how far it is from a complex line.

\smallskip
The angle between two real subspaces $V_1$ and $V_2$ of $T_y\C^2$ of the same dimension can be defined as 
\[
\ang (V_1,V_2) = \max_{u_1\in V_1} \min_{u_2\in V_2} \angle (u_1,u_2).
\]

For each $n\geq0$, let $W_n$ be the set of points from $W$ that stay in $W$ under  the first $n$ backward iterates of $f$. Let $U_n=\phi^{-1}(W_n)$.

\begin{prop}\label{prop:span} Let $x\in W_n$ and  $v\in T_{x}\Wloc$. There exists $\rho<1$ such that  
\[
\ang \left(T_{x}\Wloc, {\rm Span}_{\C}\{v\}\right) = \bigO{(\rho^n)}.
\]
\end{prop}
\proof
Let $y=f^{-n}(x)$. Let $w=df^{-n}_x(v)\in T_y \Wloc$. The vector  $J_{y} w$ does not in general belong to $T_{y}\Wloc$, but it does belong to the horizontal cone $\mathcal{C}^{h}_y$. The derivative $df^n_y$ maps this cone into a smaller cone inside $\mathcal{C}^{h}_{x}$, with an angle opening  $\bigO{(\rho^n)}$ where $\rho<1$ as in Remark \ref{rem:rho}. The vectors $v$ and 
\[
	J_x v=J_{f^n(y)}df^n_y(w) = df^n_y(J_y w)
\] 
both belong to the cone $df^n_y(\mathcal{C}^{h}_y)$. Note that ${\rm Span}_{\C}(v)={\rm Span}_{\R}(v, J_x v)$. In conclusion the angle between the complex line spanned by $v$ and $J_x v$ and the real tangent space $T_x \Wloc$ is $\bigO{(\rho^n)}$. 
\qed

Define the norm of the $\dbar_{J'}$-derivative of $f$ on a set $W$ as 
\[
\|\dbar_{J'}f\|_W = \sup_{z\in W}\|(\dbar_{J'}f)_z\|,
\] where $\|(\dbar_{J'}f)_z\|$
is the operator norm of $(\dbar_{J'}f)_z:T_zW\rightarrow T_{f(z)}W$.

\begin{lemma}\label{estimate:f} There exists a constant $C$ such that for every $n\geq 1$, 
\[
\| \dbar_{J'} f\|_{W_n} <C \rho^n.
\]
\end{lemma}
\proof Let $x\in W_n$. Then $f(x)\in \Wloc$. Let $v\in T_x\Wloc$. 
We will use the fact that the map $f$ is analytic on $B'$, so $J_{f(x)}\circ df_x  = df_x \circ  J_{x}$, to estimate
\begin{eqnarray}\label{eq:JJ}
	2\| (\dbar_{J'} f)_x (v) \| &=& \| J'_{f(x)}\circ df_x v  - df_x \circ  J'_{x}v \| \\ \nonumber
				&\leq& \| J'_{f(x)}\circ df_x v - J_{f(x)}\circ df_x v\| + \|df_x \circ  J_{x}v - df_x \circ  J'_{x}v\|.
\end{eqnarray}
Let $\beta = \angle (J_{x}v, J'_{x}v)$.  Since $J_{x}=J'_{x}$ on $\Lambda$, we may assume that $\beta<\pi/2$ for $x\in W_n$. Note then that $\beta$ is also equal to the $\ang (J_{x}v, T_{x}\Wloc)$, as $T_{x}\Wloc=\mbox{Span}_{\R}\{v, J'_{x}v\}$, and $\langle v,J'_{x}v\rangle=\langle v,J_{x}v\rangle=0$.

By a direct computation we get 
\begin{equation}\label{eq:beta1}
	\|J'_{x}v - J_{x}v\| = 2\|v\| \sin(\beta/2) \leq \beta \|v\|.
\end{equation}

Let $M = \sup \|df_x\|$, where the supremum is taken after $x\in f(\overline{B})$. Clearly $M<\infty$ since $f$ is $C^1$.  In view of Equation \eqref{eq:beta1}, we get 
\begin{equation}\label{eq:M}
\|df_x (  J_{x}v - J'_{x}v)\| \leq M  \|J_{x}v - J'_{x}v\| \leq \beta M \|v\|.
\end{equation}

By the same estimate \eqref{eq:beta1}, we have $\|J'_{f(x)}w - J_{f(x)}w\|\leq \beta' \|w\|$, where $w=df_x v$ and $\beta'$ is the angle between the vector $J_{f(x)}w$ and the tangent space $T_{f(x)}\Wloc$.  The vector $v$ is in the horizontal cone $\mathcal{C}_x^h$, so $\|w\| \leq \lambda_1^{-1}\|v\|$,  by partial hyperbolicity 
\eqref{eq:Cn}. 
Putting everything together, Equation \eqref{eq:JJ} becomes
\[
\| (\dbar_{J'} f)_x (v) \| \leq \frac{1}{2}(\lambda_1^{-1}\beta'+ M\beta)\|v\|.
\]

From Proposition \ref{prop:span} we know that $\beta=\bigO{(\rho^n)}$ and $\beta'=\bigO{(\rho^{n+1})}$, for some $\rho<1$. Thus there exists a constant $C$, independent of $n$, such that $\| (\dbar_{J'} f)_x \| < C\rho^n$, for every $x\in W_n$. 
\qed

We will now transport the estimates obtained for $f$ in Lemma \ref{estimate:f} to estimates for the  $\dbar$-derivative of $g$. Since $g=\phi^{-1}\circ f \circ \phi$ and $\phi$ is $(J',i)$-holomorphic, we get the following relation between the $\dbar$-derivatives of $f$ and $g$ computed with respect to the corresponding almost complex structures $J'$ and respectively $i$:
\begin{equation}
(\dbar_{i}g)_z=d\phi^{-1}_{\phi(g(z))}\circ (\dbar_{J'}f)_{\phi(z)} \circ d\phi_{z}, \ \ \mbox{for all}\ z\in U.
\end{equation}
Using Lemma \ref{estimate:f} and the fact that $d\phi$ and $d\phi^{-1}$ are bounded above, we obtain that there exists a constant $C'$ such that for every $n\geq 1$, 
\begin{equation}\label{cor:estimate-g}
\| \dbar_{i} g\|_{U_n} <C' \rho^n.
\end{equation}

\begin{cor}\label{cor:estimate-g2} There exists a constant $C'$ such that for every $n\geq 1$, 
\[
| \dbar g(z)| <C' \rho^n, \ \ \mbox{for all}\ z\in U_n.
\]
\end{cor}
\proof The proof follows directly from  \eqref{cor:estimate-g}. For completion, we give the details below. Let $z=x+i y$ and $g(x,y) = s(x,y) + i t(x,y)$. 
We have
\begin{eqnarray*}
dg_z\left(\partial_x\right) = \frac{\partial s}{\partial x} \frac{\partial }{\partial s} + \frac{\partial t}{\partial x} \frac{\partial }{\partial t}\ \ \ \mbox{and}\ \ \
dg_z\left(\partial_y\right) = \frac{\partial s}{\partial y} \frac{\partial }{\partial s} + \frac{\partial t}{\partial y} \frac{\partial }{\partial t}.
\end{eqnarray*}
The standard complex structure on $\C$ satisfies the relations in Equation \eqref{eq:iz}, so we compute 
\begin{eqnarray*}
\left(dg_z\circ i_z - i_{g(z)}\circ dg_z \right) \left(\partial_x \right) &=&  \left( \frac{\partial s}{\partial y}+  \frac{\partial t}{\partial x}\right)\frac{\partial }{\partial s}+ \left( \frac{\partial t}{\partial y}- \frac{\partial s}{\partial x}\right)\frac{\partial }{\partial t}, \\
\left(dg_z\circ i_z - i_{g(z)}\circ dg_z \right) \left( \partial_y \right) &=& -\left( \frac{\partial t}{\partial y}- \frac{\partial s}{\partial x}\right)\frac{\partial }{\partial s}+ \left( \frac{\partial s}{\partial y}+  \frac{\partial t}{\partial x}\right)\frac{\partial }{\partial t}.
\end{eqnarray*}

\noindent Using Equation \eqref{eq:dd}, the complex $\dbar$-derivative of $g$ is 
\[
\dbar g  = \frac{1}{2}\left(\frac{\partial s}{\partial x}-\frac{\partial t}{\partial y}+i\left(\frac{\partial t}{\partial x}+\frac{\partial s}{\partial y}\right)\right).
\]
By inequality \eqref{cor:estimate-g}, $\| (\dbar_ig)_z(\partial_x)\|$ and $\| (\dbar_ig)_z(\partial_y)\|$ are bounded above by  $C' \rho^n$
for $z\in U_n$, which implies that $|\dbar g(z)|\leq C' \rho^n $ for all $z$ in $U_n$.
\qed

\section{Quasiconformal conjugacy}\label{sec:qc}

\subsection{Preliminaries}\label{sec:prelim}
In this section we give a brief account of quasiconformal homeomorphisms. We refer to the classical text of Ahlfors \cite{Ahl} for a thorough treatment of quasiconformal maps. Let $U, V$ be two open sets of $\C$ and $\psi:U\rightarrow V$ be a homeomorphism such that $\psi$ belongs to the Sobolev space $W_{\rm loc}^{1,2}(U)$ (that is, $\psi$ has distributional first order derivatives which are locally square-integrable). Let $K\geq 1$. The map $\psi$ is called $K$-quasiconformal if 
\[
|\dbar \psi|\leq \frac{K-1}{K+1}|\partial \psi|
\]
almost everywhere. 

Quasiconformal homeomorphisms are almost everywhere differentiable. If $\psi$ is quasiconformal, then $\partial \psi\neq 0$  and $Jac(\psi)>0$ almost everywhere. 
The {\it complex dilatation} of $\psi$ at $z$ (or the {\it Beltrami coefficient} of $\psi$) is defined as
\[
\mu_{\psi} (z) = \frac{\bar{\partial} \psi (z) }{\partial \psi (z)}. 
\]
We have $\|\mu_{\psi}\|_{\infty}<1$, where $\|\cdot\|_{\infty}$ is the essential supremum. The  number
\[
K(\psi, z) = \frac{1+|\mu_\psi(z)|}{1-|\mu_\psi(z)|}
\]
is the {\it conformal distortion} of $g$ at $z$ (or the {\it dilatation} of $g$ at $z$). Clearly $\|\mu\|_{\infty}<1$ is equivalent to $K(\psi,z)<\infty$  almost everywhere.  

By Weyl's Lemma, if $\psi$ is $1$-quasiconformal, then it is conformal ({\it i.e.} if $\dbar \psi=0$ a.e. then $\psi$ is conformal). The composition of a $K_1$-quasiconformal homeomorphism with a $K_2$-quasiconformal homeomorphism is $K_1K_2$-quasi-conformal. 
The inverse of a $K$-quasiconformal homeomorphism is also $K$-quasiconformal. Using the chain rule for complex dilatations, we note that if $\mu_{\psi}=\mu_{\varphi}$ almost everywhere, then the composition $\psi\circ \varphi^{-1}$ is conformal.

Each measurable function $\mu: V\rightarrow \C$ with  $\|\mu\|_{\infty}<1$ is called a {\it Beltrami coefficient} in $V$. If $\mu:V\rightarrow \C$ is a Beltrami coefficient, then the {\it pull back} of $\mu$ under $\psi$ is a Beltrami coefficient  in $U$, denoted by $\psi^*\mu$.

\subsection{Quasiconformal conjugacy to an analytic map}\label{sec:qc}

Let $U, U'$ be open sets of $\C$ and $g:U\rightarrow U'$ be the orientation-preserving $C^1$ diffeomorphism on a neighborhood of $U$ with $g(0)=0$, constructed in Section \ref{sec:structure}. In this section, we will show how to change the complex structure on $U$, in order to make the function $g$ analytic. 

For  $n\geq0$, consider  the sets 
\begin{equation}\label{eq:Un}
U_n = \bigcap_{k=0}^n g^k(U)\ \ \mbox{and}\ \ \ U_{-n} = \bigcap_{k=0}^n g^{-k}(U).
\end{equation}
The set $U_n$ is the of points from $U$ such that the first $n$ backward iterates remain in $U$. The set $U_{-n}$ consists of the points of $U$ whose first $n$ forward iterates belong to $U$. Note that $U_0$ is equal to $U$. Moreover $U_{n+1}\subseteq U_n$ and $U_{-(n+1)}\subseteq U_{-n}$, for all $n\geq0$.  The set $X$ from Equation \eqref{eq:X} is equal to $U_{\infty}$ and $g^{-1}(X)\subseteq X$.
Note that the sets $U_n$ and $X$ have already been introduced in Section \ref{sec:structure}, as the preimages of the sets $W_n$ and $\Lambda$ from the center manifold $\Wloc$ under the parametrizing map $\phi$.

From the definitions given in Equation \eqref{eq:Un} we have the following invariance properties.
 
\begin{lemma}\label{prop:inclusions} For every $n\geq0$, we have:
\begin{itemize}
\item[a)] $g^{-n}(U_{n})= U_{-n}$
\item[b)] $g^{-1}(U_{n+1})\subset U_{n}$ and $g(U_{-(n+1)})\subset U_{-n}$
\item[c)] $g^{j}(U_{-n})=U_j\cap U_{-(n-j)}$, for $0\leq j \leq n$.
\end{itemize}
\end{lemma}
 \proof 
The set equalities
 \[g^{j}(U_{-n})=\bigcap_{k=0}^n g^{j-k}(U)= \bigcap_{k=0}^j g^{k}(U) \cap \bigcap_{k=0}^{n-j} g^{-k}(U)=U_j \cap U_{-(n-j)}\]
 can be used to prove part c). Taking $j=n$ yields part a).   The first inclusion in part b)
follows from the fact that $g(U_n)\cap U = U_{n+1}$, while the second relation is obtained from part c) by taking $j=1$.
 \qed

Let $\sigma_0$ denote the standard almost complex structure of the plane, represented by the zero Beltrami differential on $U$. 
The following lemma is a restatement of Corollary \ref{cor:estimate-g2} in the language of Beltrami coefficients. 

\begin{lemma}\label{lema:bounded1} There exist $\rho<1$ and $M$ independent of $n$ such that 
\[
\sup_{z\in U_{n}}|\mu_{g}(z)|<M\rho^{n}, \ \ \mbox{for all}\ n\geq 0.
\]
\end{lemma}

For each positive integer $n$, let $\mu_{g^{n}}:U_{-n}\rightarrow \C$ be the restriction of the pullback $(g^{n})^{*}\sigma_{0}$ to the set $U_{-n}$. This is the Beltrami coefficient of the $n$-th forward iterate $g^{n}$ on the set $U_{-n}$. 

\begin{lemma}\label{lema:bounded}
There exists a constant $\kappa<1$, such that for all integers $n\geq0$ we have
\[
|\mu_{g^{n}}(z)|<\kappa, \ \ \mbox{for all}\ z\in U_{-n}.
\]
\end{lemma}
\proof 
Let $n>0$ and $z\in U_{-n}$. Let $z_{j} = g^j(z)$ for $0\leq j\leq n$ denote the $j$-th iterate of $z$ under the map $g$.
By Lemma \ref{prop:inclusions} c), $z_{j}\in U_{j}$ for $0\leq j\leq n$.

We want to show that the dilatation $K(g^n, z)$ is bounded by a constant independent of $n$ and the choice of $z$. Recursively using the classical relation (see \cite{Ahl}) 
\[
K(f\circ g, z)\leq K(f,g(z))\, K(g,z),
\] 
we get the following estimate
\begin{equation}\label{eq:Kfn}
K(g^n, z)\leq \prod_{j=0}^{n-1} K(g, z_{j}).
\end{equation}
By definition and Lemma \ref{lema:bounded1} we have 
\[
K(g,z_{j}) = \frac{1+|\mu(z_{j})|}{1-|\mu(z_{j})|}\leq \frac{1+M\rho^j}{1-M\rho^j}.
\]
Since $\rho<1$, we can choose $j_0$ such that $M\rho^j<1/3$, for $j\geq j_0$. Then $K(g,z_{j})$ is bounded above by $1+3M\rho^j$, for $j\geq j_0$. 
The product of the first $j_0$ terms in Equation \eqref{eq:Kfn} is bounded by a constant $C$. 
This is an immediate consequence of the fact that $g$ is injective and orientation-preserving on a neighborhood of $U$, thus $\|\mu_g\|_{\infty}$ is bounded away from $1$. 
For $n> j_0$ we get 
\[
K(g^n,z)\leq C\prod_{j=j_0}^{n-1} (1+3M\rho^j).
\]
The infinite product $\prod (1+3M\rho^j)$ is convergent. 
There exists a constant $M'$ such that $K(g^n,z)<M'$ for all points $z\in U_{-n}$ and all $n\geq 0$. Thus $|\mu_{g^n}(z)|<(1-M')/(1+M')$ on $U_{-n}$. 
\qed

For $n>0$, let $\sigma_{n}:U_{n}\rightarrow \C$ be the restriction of the pullback $(g^{-n})^{*}\sigma_{0}$ to the set $U_{n}$. This is the Beltrami coefficient of the $n$-th backward iterate $g^{-n}$ on the set $U_{n}$. We write $\sigma_{n}=\mu_{g^{-n}}$ for simplicity. The map $g^n: U_{-n}\rightarrow U_n$ is bijective. In fact, the measurable function $\sigma_n$ is the push forward of $\sigma_0$ under $g^n$, also written as $\sigma_n=(g^n)_{*}\sigma_0$.

From the standard properties of Beltrami coefficients we have 
\begin{equation*}
\mu_{g^{-1}}(z) = -\left(\frac{\partial g(z)}{|\partial g(z)|}\right)^{2}\mu_{g}(g^{-1}(z)),
\end{equation*}
so $|\mu_{g^{-1}}(z)| = |\mu_{g}(g^{-1}(z))|$.  It follows that $|\mu_{g^{-n}}(z)| = |\mu_{g^{n}}(g^{-n}(z))|$, for all $z\in U_{n}$.
By Lemma \ref{lema:bounded} we get $\|\sigma_{n}\|_{\infty}<\kappa$, for all $n>0$. 

By a result of Sullivan \cite{Su}, a uniformly quasiconformal group is conjugate to a group of conformal transformations. We will give a direct proof of this property in our situation:

\begin{thm}\label{thm:conformal} 
The map $g^{-1}:U_1\rightarrow U_{-1}$ is quasiconformally conjugate to an analytic map.
\end{thm}
\proof Consider the measurable function $\mu:U\rightarrow \C$, given by 
\begin{equation}\label{eq:mu}
\mu = \left\{\begin{array}{ll}
		\sigma_{n}& \mbox{on}\ U_{n}-U_{n+1},\ \mbox{for}\ n\geq 0\\
		\sigma_{0}& \mbox{on}\ X.
	\end{array}\right.
\end{equation}
Then $\|\mu\|_{\infty}<1$ by Lemma \ref{lema:bounded} and the observation above. Thus $\mu$ is a Beltrami coefficient in $U$. Moreover, by construction, $\mu$ is $g^{-1}$ invariant, {\it i.e.} $(g^{-1})^{*}\mu = \mu$ on $U_1$. The standard almost complex structure is $g^{-1}$ invariant on $X$, since $\mu_g=0$ on $X$ by Lemma \ref{lema:bounded1}, which implies that $\mu_{g^{-1}}=0$ on $X$, by Definition \eqref{eq:mu} and the fact that $g^{-1}(X)\subset X$. The Beltrami coefficient $\mu$ is $g^{-1}$ invariant on $U_1-X$ by construction. To see this, let $n>0$ and pick any point $z\in U_n-U_{n+1}$. Then $\mu(z) = \sigma_n(z)$ on $U_n-U_{n+1}$, $\mu(z) = \sigma_{n-1}(z)$ on $U_{n-1}-U_{n}$,  and 
\[
g^{-1}(U_n-U_{n+1})\subset U_{n-1}-U_{n}
\]
by Lemma \ref{prop:inclusions} b). We have the following sequence of equalities
\[
(g^{-1})^{*}\mu(z)=(g^{-1})^{*}\sigma_{n-1}(z)=(g^{-n})^{*}\sigma_0(z)=\sigma_n(z)=\mu(z),
\]
which shows that $(g^{-1})^{*}\mu = \mu$ on $U_1-X$ as well.

 By the Measurable Riemann Mapping Theorem, there exists a quasiconformal homeomorphism  
$\psi:U\rightarrow \psi(U)\subset\C$ with complex dilatation $\mu_{\psi}$ equal to the Beltrami coefficient $\mu$. We choose $\psi$ such that $\psi(0)=0$. Let $\Omega = \psi(U_{-1})$ and $\Omega' = \psi(U_{1})$ and consider the map $h:\Omega \rightarrow \Omega'$, $h=\psi\circ g\circ \psi^{-1}$, as in the diagram below
\[
\diag{(U_1, (g^{-1})^*\mu)}{(U_{-1},\mu)}{\Omega'}{\Omega}{g^{-1}}{\psi}{\psi}{h^{-1}}
\]
\smallskip

The map $\psi\circ g^{-1}$ is a composition of two quasiconformal maps, hence quasiconformal. It has complex dilatation $\mu_{\psi\circ g^{-1}} = (g^{-1})^{*}\mu = \mu$. Since $\mu_{\psi\circ g^{-1}} = \mu_{\psi}$, the maps $h$ and  $ h^{-1} = \psi\circ g^{-1}\circ \psi^{-1}$  are conformal. 

The  conjugacy map $\psi$ is analytic on the interior of $X$, since the Beltrami coefficient $\mu$ on $X$ is equal to the standard almost complex structure $\sigma_0$.
\qed

\subsection{Proof of Theorem \ref{thm:QC}}\label{sec:proof} Let $f$ be a holomorphic germ of diffeomorphisms of $(\C^2,0)$ with a semi-indifferent fixed point at the origin. Let $\lambda= e^{2\pi i \theta}$ be the neutral eigenvalue of $0$. 
The restriction of the map $f$ to the center manifold, $f:W\rightarrow W'$ is conjugate to the map $g:U\rightarrow U'$ by a $C^1$-diffeomorphism $\phi$.  By Theorem \ref{thm:conformal}, and eventually replacing $U$ with $U_{-1}$, the map $g:U\rightarrow U'$ is  conjugate to an analytic map $h:(\Omega,0)\rightarrow (\Omega',0)$ by a quasiconformal homeomorphism $\psi$.

Let $\Lambda$ be defined as in Equation \eqref{eq:Lambda}. By Proposition \ref{horiz} , the interior of the set $\Lambda$ is a (not necessarily connected) complex submanifold of $\C^2$. To show that the conjugacy map $\psi\circ\phi^{-1}$ is analytic on the interior of $\Lambda$, we first recall from Corollary \ref{cor:submanC} that $\phi:\mbox{int}(X)\rightarrow\mbox{int}^{c}(\Lambda)$ is holomorphic, where  $X=\phi^{-1}(\Lambda)$.
From the proof of Theorem \ref{thm:conformal} it follows that map $\psi$ is analytic on the interior of the set $X$. So the composition $\psi\circ\phi^{-1}$ is holomorphic on $\mbox{int}^{c}(\Lambda)$.

The last thing we need to show is that $h'(0)=\lambda$. The map $d\phi^{-1}$ maps the complex eigenspace $E^c_0$ of the eigenvalue $\lambda$ to the complex tangent space $T_0U$. The action of $df_0$ on 
$E^c_0$ is multiplication by $\lambda$ and $dg_0=d\phi^{-1}_0\circ df_0\circ d\phi_0$, so $dg_0$ as a real matrix is just a rotation matrix of angle $\theta$. Since the tangent space $T_0U$ is complex, we can view $dg_0$ as a complex function, which means that $g'(0)=\lambda$. 

\begin{remark}
It is well-known that the multiplier of an indifferent fixed points is a quasiconformal invariant (even a topological invariant by a theorem of Naishul \cite{N}). Namely, if $g_1$ and $g_2$ are two quasiconformally conjugate (topologically conjugate) holomorphic germs of $(\C,0)$ with indifferent fixed points at the origin,  then $g_1'(0) = g_2'(0)$.  We cannot use this fact in our setting to conclude that $g'(0)=h'(0)$, because the map $g$ is not holomorphic. 
\end{remark}

To show that the two rotation numbers coincide, we use a generalization of Naishul's Theorem in 1D, due to Gambaudo, Le Calvez, and P\'ecou \cite{GLP}, which says in particular that the multiplier at the origin is a topological invariant for the class of orientation-preserving homeomorphisms of the plane which are differentiable at the origin and for which the derivative at the origin is a rotation.

\smallskip
Alternatively, to show that $h'(0)=\lambda$, we can use the hedgehog constructed in \cite{FLRT}. Let $V$ be a neighborhood of the origin in $\Wloc$ compactly contained in $W$. Let $K_f$ be a  compact, connected, and completely invariant set for $f$ given by Theorem \ref{thm:FLRT} such that $0\in K_f$ and $ K_f\cap \partial V\neq 0$.  By construction $K_f\subset V$ and $K_f$ is full in $\Wloc$ (that is, $\Wloc-K_f$ is connected). As in \cite{FLRT}, we can associate to $(f,K_f)$ an orientation-preserving homeomorphism $\tilde{f}$ of the unit circle, with rotation number $\theta$, by uniformizing the complement of $K_f$ inside the global center manifold $W^c(0)$ identified with $\R^2$, by the complement of the unit disk in $\C$. 
The set $K_h = \psi\circ\phi^{-1}(K_f)$ is a nontrivial compact, connected, full, and completely invariant set for $h$, containing the origin. Moreover, $K_h$ intersects the boundary of $V' = \psi\circ\phi^{-1}(V)$. 
By \cite[Lemma~III.3.4]{PM1},  we can associate to $(h, K_h)$ an orientation-preserving diffeomorphism $\tilde{h}$ of the unit circle, with rotation number $\theta'$, the argument of $h'(0)$. The conjugacy map $\psi\circ\phi^{-1}$ is uniformly continuous on a neighborhood of $K_f$. Therefore it defines a homeomorphism from the set of prime ends of $W^c(0)\cup\{\infty\}-K_f$ to the set of prime ends of $\hat{\C}-K_h$. We use this to obtain a conjugacy between the homeomorphisms $\tilde{f}$ and $\tilde{h}$ on the unit circle. Hence they have the same rotation number $\theta=\theta'$. 

\section{Dynamical consequences of Theorem \ref{thm:QC}}\label{sec:cons}

In this section we give the proofs of the remaining theorems stated in the introduction.

Theorem \ref{thm:B}, respectively Corollary \ref{corB1} follows from Theorem \ref{thm:PM3},  respectively Corollary \ref{cor:PM3} by applying Theorem \ref{thm:QC}. We now proceed with the proofs of Theorems \ref{thm:C}, \ref{thm:D},  \ref{thm:wandering} and Corollary \ref{cor:lin}.

\medskip
\noindent \textbf{Proof of Theorem \ref{thm:C}.} Consider a domain $B'\subset \C^{2}$ such that $B, f(B)$ are compactly contained in $B'$ and $f$ is partially hyperbolic on $B'$. Let $\Wloc$ be a center manifold of $0$ constructed with respect to the bigger set $B'$. 
Let $W$ be the connected component of $\Wloc\cap B \cap f^{-1}(B)$ containing the origin. 
By Theorem \ref{thm:QC}, the map $f|_{W}$ is quasiconformally conjugate to a holomorphic map $h:\Omega\rightarrow \Omega'$ whose multiplier at the origin is equal to the neutral eigenvalue of $f$. We can choose $B$ so that $\Omega\subset \C$ is a Jordan domain with $C^{1}$ boundary ({\it i.e.} an admissible neighborhood in P\'erez-Marco's terminology). Let $\phi:W\rightarrow \Omega$ denote the quasiconformal conjugacy. 

 Let $\mathcal{H}\subset \overline{B}$ be a hedgehog for $f$ such that $\mathcal{H}\cap \partial{B}\neq \emptyset$.  Consider a point $x\in B$, which does not belong to $W^{ss}_{\rm loc}(\mathcal{H})$. Suppose that $f^n(x)\in B$ for $n\geq 1$. By the shadowing property from Proposition \ref{prop:shadow} there exists $y\in \Wloc-\mathcal{H}$ such that $f^n(y)\in \Wloc$ for all $n\geq 1$ and the orbit of $y$ shadows the orbit of $x$. 
Clearly the $\omega$-limit set of $x$ is the same as the $\omega$-limit set of $y$. The set $K = \phi(\mathcal{H})$ is a hedgehog for $h$. The point $z=\phi(y)$ does not belong to $K$.  However, $f^{n}(z)\in \Omega$ for all $n\geq 0$. By Theorem \ref{thm:PM4}, $\omega(z)\cap K=\emptyset$. It follows that $\omega(x)\cap \mathcal{H}=\emptyset$. 

A similar argument shows that if $x\in B-\mathcal{H}$ and $f^{-n}(x)\in B$ for all $n\geq 1$, then $\alpha(x)\cap \mathcal{H}=\emptyset$.
\qed

\noindent \textbf{Proof of Theorem \ref{thm:D}.} 
Suppose that $f$ is conjugate in a neighborhood of the origin to the map $\tilde{f}(x,y)=(\lambda x, \mu(x) y)$. For some small $r$,  the disk $\D_r\times\{0\}$ is invariant under $\tilde{f}$, therefore there exists an embedded holomorphic disk $\Delta$ which is invariant for the dynamics of $f$. All local center manifolds must contain $\Delta$, therefore $0\in \Delta \subset \mathcal{H}$.

Conversely, suppose that $0\in {\rm int}^c(\mathcal{H})$, the interior of $\mathcal{H}$ relative to a center manifold. 
Let $\Delta$ be the connected component of ${\rm int}^c(\mathcal{H})$ which contains $0$. By the properties of the hedgehog $\mathcal{H}$ from Theorem \ref{thm:FLRT}, the set $\Delta$ is open, bounded, simply connected, with $f(\Delta)=\Delta$. By Proposition \ref{prop:dbar-X}, the set $\Delta$ is an analytic submanifold of $\C^2$ of dimension $1$, hence it is biholomorphic to the unit disk $\D$. Choose a biholomorphism $\phi:\D\rightarrow \Delta$ with $\phi(0)=0$, and set $g:=\phi^{-1}\circ f \circ \phi$. The map $g:\D\rightarrow \D$ is an automorphism of the unit disk, satisfying $g(0)=0$ and $g'(0)=\lambda$, $|\lambda|=1$, hence by the Schwarz Lemma we can conclude that $g$ is a rotation, $g(z)=\lambda z$ for all $z\in\D$. Therefore, the restriction of $f$ to $\Delta$ is linearizable. 

We will show that the map $f$ in conjugate in a neighborhood of the origin in $\C^2$ to a linear cocycle $(x,y)\mapsto(\lambda x, \mu(x) y)$, in a two-step argument. In the first step, we conjugate $f$ to a skew product $(x,y)\mapsto(\lambda x, \nu(x,y))$, where $\nu$ is a nonlinear cocycle, using ideas from \cite[Proposition~6]{BS}. 
In the second step we show how to reduce the nonlinear cocycle to a linear one.

\smallskip
\noindent{\bf Step 1.} We conjugate $f$ to a map of the form 
\begin{equation}\label{eq:form1}
F(x,y)=\left(\lambda x, \nu(x,y)\right),
\end{equation}
where 
\begin{equation*}\label{eq:form2}
\nu(x,y)=\mu y \left(B_{0}(x)+yB_{1}(x,y)\right),
\end{equation*}
for some holomorphic functions $B_{0}(\cdot)$ and $B_{1}(\cdot,\cdot)$ with $B_{0}(x) = 1 +\bigO(x)$. 
\smallskip

The parametrizing map $\phi:\D\rightarrow \Delta$, $\phi=(\phi_1,\phi_2)$, is a biholomorphism, hence $\left(\phi'_1(x),\phi'_2(x)\right)\neq (0,0)$ for all $x\in\D$. There exists two holomorphic maps $q_1$ and $q_2$ such that $\phi'_1(x)q_2(x)-\phi'_2(x)q_1(x)=1$ for all $x\in\D$. Let $s:\D\times\C\rightarrow\C^2$ be given by $s(x,y)=\left(\phi_1(x)+y q_1(x), \phi_2(x)+y q_2(x)\right)$, and set $F:=s^{-1}\circ f\circ s$. The maps $s$ and $F$ are local diffeomorphisms in a neighborhood of the disk $s^{-1}(\Delta)=\D\times\{0\}$, since ${\rm det}\ ds(x,0)=1$ for all $x\in \D$. Moreover $F(\D\times\{0\})=\D\times\{0\}$ and $F(x,0)=(\lambda x, 0)$ for all $x\in \D$. 

Consider the strong stable set $W^{ss}_{\rm loc}(\Delta)$ of $\Delta$ with respect to a small neighborhood $\mathcal{N}$ of $\Delta$ and let $V=s^{-1}(W^{ss}_{\rm loc}(\Delta))$. The strong stable set $W^{ss}_{\rm loc}(\Delta)$ consists of point from $\mathcal{N}$ which converge asymptotically exponentially fast to the invariant disk $\Delta$. In Step 1, we show how to straighten the foliation of the local strong stable set of $\Delta$ so that the local strong stable manifolds $W^{ss}_{\rm loc}(x)$, $x\in\Delta$, become vertical. 
 By construction, $F(V)\subset V$. The sequence of iterates $(F^{n}|_{V})_{n\geq 0}$ is a normal family. Therefore there exists a subsequence of iterates $F^{n_j}$ which converges uniformly on compact subsets of $V$ to a holomorphic map $\rho:V\rightarrow \D\times\{0\}$. One can in fact choose the subsequence $n_j$ as in \cite{B} so that the map $\rho$ is a retract of $V$ onto the invariant disk $\D\times\{0\}$, that is, $\rho(x,0)=(x,0)$ for all $x\in\D$. By construction, $\rho$ commutes with the map $F$, so we have 
\[
\rho \circ F(x,y)= F\circ \rho (x,y)=\lambda \rho(x,y).
\]

Consider the map $H(x,y)=(\rho(x,y), y)$ which leaves the disk $\D\times\{0\}$ invariant and is invertible in a neighborhood of this disk since $\partial_x \rho (x,0)=1$. By replacing $F$ with $H\circ F \circ H^{-1}$, we may assume that $F$ is linear in the first coordinate and therefore has the form given in Equation \eqref{eq:form1}. 

\smallskip
\noindent{\bf Step 2.} Next we show that $F$ is conjugate to $\tilde{f}(x,y)=(\lambda x, \mu(x) y)$, where $\mu(x) = \mu + \bigO(x)$ is  a holomorphic function. 
\smallskip

Let $F^n(x,y) = (\lambda^nx, \nu_n(x,y))$ denote the $n$-th iterate of $F$, where $\nu_n$ is holomorphic and $\nu_n = \nu_{n-1}\circ F$, for all $n\geq 1$. By convention, $\nu_{0}=\nu$. Using this recurrence relation, we find 
\begin{equation}\label{eq:nu}
\nu_{n}(x,y) =  \mu^{n}y \prod_{j=0}^{n-1} \left( B_{0}(\lambda^{j}x)+\nu_{j}(x,y)B_{1}(F^j (x,y) ) \right),
\end{equation}
for all $n\geq 1$. 

Note that partial derivative of $\nu$ with respect to $y$ has the form 
\[
\partial_{y}\nu(x,y) = \mu\left(B_{0}(x) + y C_{1}(x,y)\right),
\]
for some appropriate holomorphic function $C_{1}$. By direct computation, we obtain
\[
\partial_y \nu_n  = (\partial_y \nu_{n-1}\circ F)\cdot \partial_y \nu,
\]
which yields 
\[
\partial_y \nu_n (x,y)= \prod_{j=0}^{n-1} \partial_y \nu(F^j(x,y)) = \mu^n \prod_{j=0}^{n-1} \left( B_{0}(\lambda^{j}x)+\nu_{j}(x,y)C_{1}(F^j (x,y) ) \right),
\]
for all $n\geq 1$. Note that $\nu_{j}(x,0) = 0$ for all $j\geq 0$. Also $\mu\neq 0$ since $F$ is a local diffeomorphism. Hence, when $y=0$, the formula above simplifies to $\partial_y \nu_n (x,0) = \mu^n \prod_{j=0}^{n-1} B_{0}(\lambda^{j}x)$.  
We will show that the infinite product
\begin{equation}\label{eq:limit-psi}
\psi (x,y) = \lim\limits_{n\rightarrow \infty} \frac{\nu_{n}(x,y)}{\partial_y \nu_n (x,0)} = y \prod_{j=0}^{\infty}\left(1+\frac{\nu_{j}(x,y)B_{1}(F^j (x,y) )}{B_{0}(\lambda^{j}x)}\right)
\end{equation}
is uniformly convergent in some neighborhood $V\subset \C^{2}$ of $0$. Using the local dynamics, we can choose a sufficiently small neighborhood $V$ of the origin so that $F(x,y)\in V$ whenever $(x,y)\in V$. There exists a constant $M>0$ such that $|B_{0}(x)|<1+M|x|$ and $|B_{1}(x,y)|<M$ throughout $V$. Since $|\mu|<1$, we can choose $V$ small enough so that the following technical condition holds:
\begin{equation}\label{eq:condition}
|\mu|^{1/2}\left(1+|x|M+|y|M\right)<1,
\end{equation}
for all $(x,y)\in V$.

\begin{lemma}\label{lemma:bound} 
$|\nu_{n}(x,y)|<|\mu|^{n/2}|y|$ for all $n\geq 0$ and for all $(x,y)\in V$.  
\end{lemma}

\proof We proceed by induction. From the definition of $\nu$ and assumption \eqref{eq:condition}, for $n=0$ we get 
\[
|\nu_{0}(x,y)|\leq |\mu||y|(1+|x|M+|y|M)<|\mu|^{1/2}|y|. 
\]
Let $n\geq 1$ and suppose that $|\nu_{j}(x,y)|<|\mu|^{j/2}|y|$ for all $0\leq j\leq n-1$.  By Equation \eqref{eq:nu} and the fact that $|\lambda|=1$ and $|\mu|<1$ we have
\[
|\nu_{n}(x,y)|\leq |y||\mu|^{n/2}\prod\limits_{j=0}^{n-1}|\mu|^{1/2}\left(1+|x|M+|y|M\right) <|y||\mu|^{n/2},
\]
which concludes the proof.
\qed

The germ $F$ is a local diffeomorphism, so the Jacobian is bounded away from $0$. Thus there exists a constant $\kappa>0$ such that $|B_{0}(x)|>\kappa$ for all $x\in\Delta$. 
Using Lemma \ref{lemma:bound}, we find that the infinite product \eqref{eq:limit-psi} is bounded above by 
\[
|y| \prod_{j=0}^{\infty}\left(1+|\nu_{j}(x,y)|M\kappa^{-1}\right)< |y| \prod_{j=0}^{\infty}\left(1+|\mu|^{j/2}M\kappa^{-1}\right)<\infty.
\]
This shows that the product from Equation \eqref{eq:limit-psi} is convergent, uniformly on $U$. 

From the definition of $\nu_{n}$ we have 
\[
\frac{\nu_{n+1}(x,y)}{\partial_y \nu_{n+1} (x,0)} = \frac{\nu_{n}(\lambda x,\nu(x,y))}{\partial_y \nu_n (\lambda x,0)\cdot \partial_y \nu ( x,0)}.
\]
By letting $n\rightarrow \infty$ we see that the map $\psi$ satisfies the equation
\[
\psi(F(x,y)) = \mu B_{0}(x) \psi(x,y).
\] 
Let $\Psi(x,y) = (x, \psi(x,y))$ and $\tilde{f}(x,y) = (\lambda x, \mu B_{0}(x)y)$. 
The map $\Psi$ is a holomorphic function on a neighborhood of the origin with $\Psi(0,0) = (0,0)$, which conjugates $F$ to $\tilde{f}$.
For simplicity, we denote $\mu B_{0}(x)$ by $\mu(x)$. 

This step concludes the proof of Theorem \ref{thm:D}. 
\qed

Assume as in Theorem \ref{thm:D} that the germ $f$ is analytically conjugate to $\tilde{f}(x,y)=(\lambda x, \mu(x)y)$, for some holomorphic function $\mu(x)=\mu+\bigO(x)$. Since we work in the dissipative setting, we can topologically linearize $f$ in a neighborhood of the origin. Let us now discuss the analytic linearizability of the germ $f$. The map $\tilde{f}$ can be viewed as a linear cocycle.

As in \cite{KK}, we say that a cocycle $h$ is {\it reducible} if it cohomologous to a constant map, that is, if $h$ satisfies the cohomology equation
\begin{equation}\label{eq:cohom}
h(x)-h(0) = \phi(\lambda x)- \phi(x),
\end{equation}
for some function $\phi$. 
The reducibility of $h$ depends on finer arithmetic properties of the neutral eigenvalue $\lambda$. 

Suppose that there exists a biholomorphic map $(x,y)\mapsto(x,\eta(x)y)$ with $\eta(0)=1$ which conjugates $\tilde{f}$ to the linear map $(x,y)\mapsto (\lambda x,\mu y)$. 
The maps $\mu(x)$ and $\eta(x)$ are not identically vanishing in a neighborhood of the origin, so they have holomorphic logarithms $h(x)=\log (\mu(x))$ and $\phi(x)=\log(\eta(x))$ which must satisfy the cohomology equation \eqref{eq:cohom}. If we compare the Taylor series expansions of $h(x)=\log(\mu)+a_1 x + a_2 x^2 +\ldots$ and $\phi(x)=x+b_1 x +b_2 x^2+\ldots\ $, we get $a_n =  b_n (\lambda^n -1)$ for $n\geq 1$. When $\lambda$ is not a root of unity, we can solve for $b_n$ and obtain a formal series defining $\phi$. 
 
The problem of convergence of the formal series for $\phi$ is strongly related to how fast $\lambda^n-1$ approaches $0$. Let $\lambda=e^{2\pi i \alpha}$, $\alpha\notin\Q$ and let $p_{n}/q_{n}$ be the convergents of $\alpha$ given by the continued fraction. If $\lambda$ satisfies the Brjuno condition,
\begin{equation}\label{eq:Brjuno}
\sum\limits_{n\geq 0}\frac{\log q_{n+1}}{q_{n}} < \infty,
\end{equation} 
then by \cite{Brj} we have a convergent power series for $\phi$. On the other hand, for general functions $\mu(x)$ and neutral eigenvalues $\lambda$ which are not roots of unity and do not satisfy the Brjuno condition,  small divisor problems could prevent the existence of solutions of the cohomology equation \eqref{eq:cohom}. 

\medskip
\noindent{\bf Non-linearizable germ with a Siegel disk.} 
Let $\alpha$ be an irrational angle such that 
\[
\limsup\limits_{n\rightarrow\infty}(\{n\alpha\})^{-1/n}=\infty,
\] 
where $\{n\alpha\}$ denotes the fractional part of $n\alpha$. Set $\lambda=e^{2\pi i \alpha}$. As in \cite{M}, the arithmetic condition imposed on $\alpha$ is equivalent to
\[
\limsup\limits_{n\rightarrow\infty}|\lambda^n-1|^{-1/n}=\infty.
\] 
Assume that the power series expansion of $h(x)=\log(\mu(x))$ has radius of convergence $0<R<\infty$. Since $b_n=a_n(\lambda^n-1)^{-1}$ for $n\geq 1$, it follows that the radius of convergence of the formal series expansion of $\phi$ is $0$,
which shows that there is no holomorphic function $\phi$ satisfying  \eqref{eq:cohom}. 
As a concrete example, let $\lambda$ be chosen as above, and let $\mu(x)=\mu e^{x/(1-x)}$, with $|\mu|<1$. 
Then $\tilde{f}(x,y)=(\lambda x, \mu(x)y)$ is an example of a local diffeomorphism which has a Siegel disk containing $0$ and which cannot be linearized in a neighborhood of the origin in $\C^2$.

\medskip
The next two proofs in this section deal with the case of polynomial automorphisms of $\C^2$.

\medskip
\noindent \textbf{Proof of Corollary \ref{cor:lin}.} 
If $f$ is a polynomial automorphism of $(\C^2,0)$ then it has non-zero constant Jacobian, equal to the product $\lambda \mu$ of the two eigenvalues of $df_{0}$. By Theorem \ref{thm:D}, $0\in {\rm int}^c(\mathcal{H})$ if and only if $f$ is analytically conjugate to a holomorphic map $\tilde{f}(x,y)=(\lambda x,\mu(x)y)$. Any conjugacy function constructed in the proof of Theorem \ref{thm:D} has the property that the determinant of its Jacobian matrix restricted to the invariant disk $\D\times\{0\}$ is $1$. This is obvious for the conjugacy maps considered at Step 1. For the coordinate transformation $\Psi$ from Step 2, this property follows from the fact that 
\[
{\rm det}\ d\Psi(x,y) = \partial_y\psi(x,y)= \lim\limits_{n\rightarrow \infty} \displaystyle \frac{\partial_y\nu_{n}(x,y)}{\partial_y \nu_n (x,0)},
\]
 which is equal to $1$ when $y=0$. Therefore ${\rm det}\ d\tilde{f}|_{\D\times\{0\}}$ is constant and equal to $\lambda\mu$. It follows that $\mu(x)=\mu$ for all $x\in \D$, so $\tilde{f}$ is a linear map.
 \qed
 
\noindent \textbf{Proof of Theorem \ref{thm:wandering}.} 
 Let $f$ be a  polynomial diffeomorphism of $\C^2$ with an irrationally semi-indifferent fixed point at the origin. Since $f$ is assumed non-linearizable in a neighborhood of the origin in $\C^2$, by Corollary \ref{cor:lin}, we know that $0\notin {\rm int}^{c}(\mathcal{H})$, the interior of $\mathcal{H}$ relative to a center manifold $\Wloc$. 
By Theorem \ref{thm:QC}, the restriction of $f$ to $\Wloc$ is quasiconformally conjugate to a holomorphic diffeomorphism $h:(\Omega,0)\rightarrow (\Omega',0)$, $h(z) = \lambda z + \bigO(z^2)$. Denote by $\phi:\Wloc\rightarrow \Omega$ the conjugacy map and by $K=\phi(\mathcal{H})$ the corresponding hedgehog for $h$. It follows that $0\notin {\rm int} (K)$, which is equivalent, by Theorem \ref{thm:PM1}, to the fact that $h$ is non-linearizable as well. By \cite{PM3}, the interior of a non-linearizable hedgehog is empty, therefore ${\rm int}^c(\mathcal{H})=\emptyset$. Hence $\mathcal{H}$ belongs to the Julia set $J$. 

We now show that a stronger statement holds true, {\it i.e.} $\mathcal{H}\subset J^{*}$. The following argument is due to Romain Dujardin. To fix the notations, let $B\Subset B'$ be two nested domains in $\C^2$ containing the fixed point $0$ such that $f$ is partially hyperbolic in a neighborhood $B''$ of $\overline{B'}$. We can choose $B'$ as in Theorem \ref{thm:C}. Let $\mathcal{H}$ and $\mathcal{H'}$ be the hedgehogs constructed with respect to $B$ and $B'$. Then $\mathcal{H}\subset \mathcal{H'}$. 

Let $x\in \mathcal{H}$. The local strong stable manifold $W^{ss}_{\rm loc}(x)$ is a holomorphic disk. The global strong stable manifold $W^{ss}(x) = \bigcup_{n\geq 0}f^{-n}(W^{ss}_{\rm loc}(f^{n}(x))$ is biholomorphic to $\C$ (see \cite[Proposition~7.3]{HOV}). Let $p$ be a saddle periodic point of $f$. By \cite{BLS}, there exists a transverse intersection point $t$ between $W^{ss}(x)$ and $W^{u}( p )$. Then $t\in f^{-k}(W^{ss}_{\rm loc}(f^{k}(x))$ for some $k\geq 0$. Iterating forward we obtain $f^{n}(t) \in W^{ss}_{\rm loc}(f^{n}(x))$ for every $n\geq k$. Using the recurrence of the points of the hedgehog from Corollary \ref{corB1}
and the strong  contraction along the leaves of $W^{ss}_{\rm loc}(\mathcal{H})$, we conclude that there exist transverse intersection points $z$ arbitrarily close to $x$. Therefore, to show that $x\in J^{*}$, it suffices to show that $z\in J^{*}$. 

We claim that $z$ belongs to the relative boundary of $K^{+}$ in $W^{u}( p )$, thus $z\in J^{*}$ by \cite[Lemma~5.1]{DL}. Otherwise, assume by contradiction that there exists a holomorphic disk $\Delta\subset W^{u}( p ) \cap K^{+}$ which contains $z$. The sequence of iterates $(f^{n}|_{\Delta})_{n\geq0}$ is a normal family.
By reducing the size of $\Delta$, we may assume that all forward iterates of $\Delta$ are contained in $B'$. 
Consider the limit $\Gamma$ of a convergent subsequence $(f^{n_{j}}(\Delta))_{j\geq 0}$. 
For any fixed $k\in \Z$,
$f^{k}(\Gamma) = \lim\limits_{j\rightarrow\infty} f^{n_{j}+k}(\Delta)\subset B'$, so $\Gamma$ belongs to every center manifold $\Wloc$ defined relative to $B''$. By Theorems \ref{thm:QC} and \ref{thm:PM2} it follows that $\Gamma\subset \mathcal{H'}$.
Therefore $\Delta\subset W^{s}_{\rm loc}(\mathcal{H'})=W^{ss}_{\rm loc}(\mathcal{H'})$ by Corollary \ref{corC2}. On the other hand, notice that $\Delta-W^{ss}_{\rm loc}(\mathcal{H'})$ is nonempty, since the holonomy is a local homeomorphism and the hedgehog $\mathcal{H'}$ has empty relative interior. This is a contradiction.

To show the last part of the theorem, suppose there is a wandering component converging to $\mathcal{H}$, and choose any interior point $z$ of the wandering component. Then for $n$ sufficiently large, we may assume that all points $f^n(z)\in B'$ and  $\omega(z)\subset\mathcal{H}'$. This again contradicts Theorem \ref{thm:C}.
\qed

Suppose $f$ is a  polynomial diffeomorphism of $\C^2$ with a semi-Siegel fixed point at $0$ with an eigenvalue $\lambda=e^{2\pi i \alpha}$, $\alpha\notin \Q$. 
Let $\varphi:\D\rightarrow \C^2$ be a  maximal injective holomorphic mapping such that $f(\varphi(\xi))=\varphi(\lambda \xi)$, for all $\xi\in\D$. The image $\Delta = \varphi(\D)$ is called a {\it Siegel disk}.

\begin{prop}\label{prop:Siegel} 
If the closure of the Siegel disk $\Delta$ is contained in a hedgehog $\mathcal{H}$, then $\partial\Delta\subset J^{*}$. 
\end{prop}
\proof 
The first observation is that since the closure of the Siegel disk is contained in a region where the map $f$ is partially hyperbolic, every point $x$ on the boundary of the Siegel disk has a well-defined local strong stable manifold, $W^{ss}_{\rm loc}(x)$ which is a holomorphic disk. Let $x\in \partial \Delta$ and let $p$ be a saddle periodic point of $f$. 
Using the same arguments as in the proof of Theorem \ref{thm:wandering}, we deduce that there exist transverse intersection points $z$ between $W^{ss}_{\rm loc}(f^{n}(x))$ and $W^{u}( p)$, arbitrarily close to $x$.  To show that $x\in J^{*}$, it suffices to show that $z$ belongs to the relative boundary of $K^{+}$ in $W^{u}( p )$, thus $z\in J^{*}$. 

Otherwise, there exists a holomorphic disk $D\subset W^{u}( p ) \cap K^{+}$ which contains $z$.
The sequence of iterates $(f^{n}|_{D})_{n\geq0}$ is a normal family, so
by reducing the size of $D$, we may assume that all its forward iterates are contained in $B'$. 
Let $g$ be the limit of some convergent subsequence $(f^{n_{j}}|_{D})_{j\geq 0}$. Denote by $\Gamma=g(D)$ and let $x'=g(z)\in \partial\Delta$.  The rank of $g$ can be $1$ or $0$, and in both cases the set $\Gamma$ belongs to the hedgehog $\mathcal{H'}$ and so to every local center manifold defined relative to $B''$. If the rank of $g$ is 1 then $\Gamma$ is a holomorphic disk which extends the Siegel disk to a neighborhood of $x'$ in $\Wloc$; this is a contradiction to the maximality of the Siegel disk. If the rank of $g$ is 0 then $\Gamma$ is a single point $x'$. However, this is not possible, since $D\cap W^{ss}_{\rm loc}(\Delta)\neq \emptyset $ by construction, and if we pick any point in the intersection, then its limit under the subsequence $f^{n_{j}}$ will be an interior point of the Siegel disk, and not the boundary point $x'$. Hence in both cases we reach a contradiction.
\qed

We proceed with the proof of Theorem \ref{thm:F}, which discusses germs with a semi-parabolic fixed point at the origin. 
\medskip

\noindent \textbf{Proof of Theorem \ref{thm:F}.} Let $B\subset \C^2$ be a ball containing $0$ such that $f$ is partially hyperbolic on a neighborhood $B'$ of $\overline{B}$.
Let $\Wloc$ be any local center manifold of $0$, constructed relative to $B'$. By the weak uniqueness property of center manifolds, $\Sigma_{B}\subset \Wloc$, where $\Sigma_B$ is defined in Equation \eqref{eq:Sigma}. By Theorem \ref{thm:QC}, the map $f$ restricted to $\Wloc$ is quasiconformally conjugate to an analytic map $h:(\Omega,0)\rightarrow(\Omega',0)$, where $\Omega,\Omega'$ are domains in $\C$. By Corollary \ref{cor:submanC}, the quasiconformal map $\phi^{-1}$ is holomorphic on the interior of $\Lambda$ rel $\Wloc$, for the set $\Lambda$ defined in \eqref{eq:Lambda}. Therefore $\phi^{-1}$ is holomorphic on the set $\Sigma_B$, since $\Sigma_B$ belongs to the interior of $\Lambda$ rel $\Wloc$.

By Theorem \ref{thm:QC}, the map $h$ has a parabolic fixed point at $0$, with multiplier $\lambda$, so it is conjugate to a normal form 
\[
h(z)=\lambda z + z^{\nu q +1}+a z^{2\nu q +1}+\bigO(z^{2\nu q+2}).
\] 

By the Leau-Fatou theory of parabolic holomorphic germs of $(\C,0)$, $h$ has $\nu$ cycles of $q$ attracting and $q$ repelling petals, containing $0$ in their boundary. On each repelling petal $\mathcal{P}_{\rm rep}$, there exists an outgoing Fatou coordinate $\varphi^{o}:\mathcal{P}_{\rm rep}\rightarrow \C$, which satisfies the Abel equation $\varphi^{o}(f^{q})=\varphi^{o}+1$, that is, it conjugates $f^{q}$ to the translation $z\mapsto z+1$ on a left half plane.

The repelling petals for $f$ are just the pull-back of the repelling petals for $h$ under the holomorphic map $\phi^{-1} |_{\Sigma_B}$, and on each such repelling petal we have a holomorphic Fatou coordinate $\varphi^{o}\circ\phi^{-1}$.
\qed

We define the parabolic basin of $0$ with respect to the neighborhood $B$ as
\[
\mathcal{B}_{par}(0) = \{x\in B : f^{n}(x)\in B\ \forall n\in \N, f^{n}(x)\rightarrow 0 \mbox{ locally uniformly}\}.
\]
Note that this set has complex dimension two. Therefore we cannot use the same strategy as in the proof of Theorem \ref{thm:F} to construct incoming Fatou coordinates, since the conjugacy map $\phi^{-1}$ is only quasiconformal on the set one-dimensional slice $\Wloc \cap \mathcal{B}_{par}(0)$.

\medskip
We conclude this section with a generalization of Naishul's theorem to higher dimensions, which is of independent interest. 

\begin{thm}\label{thm:Naishul}  For $i=1,2$, let $f_i$ be a holomorphic germ of diffeomorphisms of $(\C^2,0)$ with a fixed point with eigenvalues $\lambda_i=e^{2\pi i \theta_i}$ and $|\mu_i|<1$. If $f_1$ and $f_2$ are topologically conjugate by a homeomorphism $\varphi:\C^{2}\rightarrow \C^{2}$ with $\varphi(0)=0$, then $\theta_{1}=\pm \theta_{2}$.   
\end{thm}
\proof 
Let $\varphi$ be a homeomorphism which fixes the origin and conjugates $f_{1}$ and $f_{2}$, that is $\varphi \circ f_1 = f_2\circ \varphi$. Let $B_1'\subset \C^{2}$ be a small enough ball containing the origin such that $f_{1}$ is partially hyperbolic on a neighborhood of $\overline{B_1'}$ and $f_{2}$ is partially hyperbolic on a neighborhood of the closure of $B'_2=\varphi(B_1')$. Consider a  ball $B_{1}\Subset B_1'$ and let $\mathcal{H}_{1}$ be a hedgehog for $f_{1}$ and $B_{1}$. Let $\mathcal{H}_{2}=\varphi(\mathcal{H}_{1})$ be the corresponding hedgehog for $f_{2}$ and $B_{2}=\varphi(B_{1})$.

Let $i=1$ or $2$.  Fix $r\geq 1$ and let $W^{c}_{i}$ be a $C^{r}$-smooth local  center manifold of $0$ for $f_{i}$, constructed relative to $B'_{i}$. 
The hedgehog $\mathcal{H}_{i}$ is compactly contained in the center manifold $W^{c}_{i}$.  By Theorem \ref{thm:FLRT}, there exists a dynamically defined strong stable lamination  $W^{ss}_{\rm loc} ( \mathcal{H}_i )$ invariant under $f_i$, whose leaves are holomorphic disks. Note that the conjugating homeomorphism $\varphi$ preserves the strong stable lamination of the hedgehogs.

By the general theory of partially hyperbolic systems (see e.g. \cite{SSTC} pp. 278-282), there exists a $C^{r-1}$-smooth invariant foliation $\mathcal{F}^{ss}_{\rm loc}$ 
with $C^r$-smooth leaves transverse to the center manifold $W^{c}_{i}$, which extends the strong stable lamination $W^{ss}_{\rm loc} ( \mathcal{H}_i )$. 
Consider the map $\pi: B_{2}\rightarrow W^{c}_{2}$ obtained by projecting along the leaves of the stable foliation to $W^{c}_{2}$. The projection $\pi$ is continuous. There exist neighborhoods $U_{i}$ of $\mathcal{H}_{i}$ in $W^{c}_{i}$ such that the  homeomorphism  $\varphi$ induces a continuous surjection $\varphi^{c}:U_{1}\rightarrow U_{2}$ given by $\varphi^{c}=\pi\circ \varphi$. The map $\varphi^{c}$ is a homeomorphism from $\mathcal{H}_{1}$ to $\mathcal{H}_{2}$ which conjugates $f_1|_{\mathcal{H}_{1}}$ to $f_2|_{\mathcal{H}_{2}}$, but it might not conjugate $f_1|_{U_{1}-\mathcal{H}_{1}}$ to $f_2|_{U_{2}-\mathcal{H}_{2}}$.

By Theorem \ref{thm:QC} there exists a quasiconformal map $\phi_{i}:W_{i}^{c}\rightarrow \Omega_{i}\subset\C$ which conjugates $f_{i}|_{W_{i}^{c}}$ to an analytic map $h_i:\Omega_{i}\rightarrow \Omega_{i}'$, with multiplier $h_i'(0)=e^{2\pi i \theta_{i}}$. Let $K_{i}=\phi_{i}(\mathcal{H}_{i})$ be the corresponding hedgehog for the map $h_{i}$ and the Jordan domain $D_i=\phi_{i}(B_{i}\cap W^{c}_{i})$. 
By the Uniformization Theorem, there exist a Riemann map $\psi_{i}:\hat{\C}-\overline{\D}\rightarrow \hat{\C}-K_{i}$ which fixes $\infty$. When $K_i$ is not locally connected, $\psi_{i}$ does not extend as a continuous map from $\s^1$ to $K_i$, but has radial limits almost everywhere on $\s^1$.
Let $g_{i} = \psi_{i}^{-1}\circ h_{i}\circ\psi_{i}$ 
be the corresponding holomorphic map induced by $f_{i}$ in an annular neighborhood of $\s^{1}$ in $\hat{\C}-\overline{\D}$. 
By \cite[Theorem~2]{PM1}, the map $g_{i}$ extends to an analytic circle diffeomorphism with rotation 
number  $\theta_{i}$.

The map $\eta = \phi_{2}\circ \varphi^{c}\circ \phi_{1}^{-1}:\Omega_{1}\rightarrow \Omega_2$ is continuous. It is a homeomorphism from $K_{1}$ to $K_{2}$ conjugating $h_1|_{K_1}$ to $h_2|_{K_2}$ and it maps the complement of $K_1$ to the complement of $K_2$. The following result follows from the theory of prime ends:

\begin{lemma}\label{prime-ends} Let  $\eta: K_1 \to K_2$ be a homeomorphism between full compact sets $K_i\subset\C$ which extends continuously to a neighborhood $\Omega_1\supset K_1$ such that  $\eta (\Omega_1- K_1) \subset  (\C- K_2)$.
Then $\eta$ induces a continuous map $\varrho: \s^1\to \s^1$ between the corresponding spaces of prime ends.
\end{lemma}

By Lemma \ref{prime-ends}, the map $\varrho= \psi_{2}^{-1}\circ \eta \circ\psi_{1}$ extends continuously to $\s^{1}$. 
We will show that $\varrho$ is a homeomorphism which conjugates $g_1$ and $g_2$ on $\s^1$, {\it i.e.} $\varrho\circ g_1 = g_2\circ \varrho$, which will prove that $\theta_1=\pm \theta_2$.

It is easy to check that $\theta_1$ is rational iff $\theta_2$ is rational. The rational case is easy to deal with and is left to the reader. 
Suppose $\theta_{1}$ is irrational. A point $z\in \partial K_{1}$ is accessible if there is a path $\gamma$ in $\C-K_{1}$ landing at $z$. 
By a theorem of Lindel\"of, $z$ is accessible if and only if it is the landing point of some ray $R_{t} = \{\psi_{1}(r e^{2\pi i t}), r>1\}$. Let $z_{1}\in K_{1}\cap\partial D_{1}$. There is a path $\gamma_{1}$ lying outside $D_{1}$ landing at $z_{1}$, so $z_{1}$ is accessible.  The image $\eta(\gamma_{1})$ is a curve in $\C-K_{2}$, so the point $z_{2} = \eta(z_{1})$ of $K_{2}\cap \partial D_{2}$ is also accessible.  
 Let $t_{1}, t_{2} \in\s^{1}$ be the angles corresponding to $z_{1}$, respectively $z_{2}$.  Then $\varrho(t_{1}) = t_{2}$ and  $\varrho$ conjugates $g_1$ to $g_2$ on the orbit of $t_{1}$ under $g_{1}$.   By Denjoy's theorem, $g_1$ is topologically conjugate to an irrational rotation, so the orbit of $t_1$ is dense in the circle. Similarly, the orbit of $t_2$ is  dense in $\s^1$, so $\varrho:(\s^1,t_1)\rightarrow (\s^1,t_2)$ semi-conjugates $g_1$ and $g_2$ on $\s^1$. By replacing $f_1$ with $f_2$ and reversing the entire construction, one obtains a semi-conjugacy $\varrho':(\s^1,t_2)\rightarrow (\s^1,t_1)$ between the two analytic circle diffeomorphisms $g_1$ and $g_2$.  It follows that $\varrho$ and $\varrho'$ must be mutually inverse conjugacies. 
\qed

We remark that Theorem \ref{thm:Naishul} generalizes in a straightforward way to the case of holomorphic germs of diffeomorphisms of $(\C^n,0)$, for $n>2$, which have a fixed point at the origin with exactly one eigenvalue on the unit circle and $n-1$ eigenvalues inside the unit disk.

\section{A generalization for germs of $(\C^n,0)$}\label{sec:Cn}

In this section we consider holomorphic germs $f$ of diffeomorphisms of $(\C^{n},0)$ such that the linear part of $f$ at $0$ has eigenvalues $\lambda_{i}$, $1\leq i\leq n$, with $|\lambda_k|=1$ and
\begin{equation}\label{eq:k}
0<|\lambda_{1}|\leq \ldots \leq |\lambda_{k-1}|<1<|\lambda_{k+1}|\leq \ldots\leq |\lambda_{n}|,
\end{equation}
for some $k$ between $1$ and $n$.

The presence of the neutral eigenvalue permits the existence of a rich type of local invariant sets and induces more complicated local dynamics.

The tangent space at $0$ has an invariant splitting into three subspaces $T_0\C^n=E_0^s\oplus E_0^c\oplus E_0^u$, of dimensions $k-1$, $1$, and respectively $n-k$. $E^s_0$ is strongly contracted and $E_0^u$ is strongly expanded by $df$, and the center direction $E_0^c$ is the eigenspace corresponding to the neutral eigenvalue $\lambda_k$. 
When $k\neq 1$ and $k\neq n$, $f$ is partially hyperbolic (in the narrow sense) (see \cite{HP}, \cite{CP}). Partial hyperbolicity is an open condition which can be extended to a suitable neighborhood of the origin. 

Let $B'$ a small ball containing the origin. As in Section \ref{sec:structure} we explain in terms of invariant cone fields what it means for $f$ to be partially hyperbolic on $B'$. Let $E$ be a subspace of $T_x\C^n$ and denote by 
\[C_x(E,\alpha)=\{ v\in T_x \C^n, \angle(v,E)\leq \alpha \}\] the cone at $x$ of angle $\alpha$ centered around $E$. 

There exist (not necessarily invariant) continuous distributions $E^s$, $E^c$ and $E^u$, extending $E_0^s$, $E_0^c$, $E_0^u$, such that $T_x\C^n=E^s_x\oplus E^c_x\oplus E^u_x$ for any $x$ in $B'$. Let $
E^{cs}_x= E^s_x\oplus E^c_x \ \ \mbox{and} \ \ \ E^{cu}_x=E^c_x\oplus E^u_x
$. There exist invariant cone families of stable and unstable cones
\begin{equation*}
\mathcal{C}_x^{s}=C_x(E^s_x,\alpha), \ \ \  \mathcal{C}_x^{u}=C_x(E^u_x,\alpha)
\end{equation*}
and center-stable and center-unstable cones
\begin{equation*}
\mathcal{C}_x^{cs}=C_x(E^{cs}_x,\alpha), \ \ \  \mathcal{C}_x^{cu}=C_x(E^{cu}_x,\alpha)
\end{equation*}
such that 
\begin{eqnarray*}
d_xf^{-1}(\mathcal{C}_x^{s})\subset \mbox{Int}\ C_{f^{-1}(x)}^{s}\cup \{0\}, && d_xf(\mathcal{C}_x^{u})\subset \mbox{Int}\ C_{f(x)}^{u}\cup\{0\}\\
d_xf^{-1}(\mathcal{C}_x^{cs})\subset \mbox{Int}\ C_{f^{-1}(x)}^{cs}\cup\{0\}, && d_xf(\mathcal{C}_x^{cu})\subset \mbox{Int}\ C_{f(x)}^{cu}\cup\{0\}
\end{eqnarray*}
and there are constants $0<\mu_s<\mu_{cu}\leq1\leq\mu_{cs}<\mu_{u}$ such that $\mu_{cu}<\mu_{cs}$, $|\lambda_{k-1}|<\mu_s$, $|\lambda_{k+1}|>\mu_u$, and the following inequalities hold:
\begin{eqnarray*}
\|df_x(v)\|  &\leq& \mu_s\,\|v\|,\ \ \ \ \mbox{for}\ v\in \mathcal{C}_{x}^{s} \\
\|df_x(v)\|  &\leq& \mu_{cs}\,\|v\|,\ \ \ \mbox{for}\ v\in \mathcal{C}_{x}^{cs} \\
 \|df_x(v)\|  &\geq& \mu_u\,\|v\|, \ \ \ \ \mbox{for}\ v\in \mathcal{C}_{x}^{u}\\
  \|df_x(v)\|  &\geq& \mu_{cu}\,\|v\|, \ \ \ \mbox{for}\ v\in \mathcal{C}_{x}^{cu}
\end{eqnarray*}
The fact that $f$ is partially hyperbolic on the set $B'$ implies that there exists local center-stable manifolds $W^{cs}_{\rm loc}$ and center-unstable manifolds $W^{cu}_{\rm loc}$  of class $C^1$, tangent at $0$ to the subspaces $E^{cs}_0$ and respectively to $E^{cu}_{0}$. By intersecting the local center-stable and the local center-unstable manifolds, one shows the existence of center manifolds $\Wloc$ of class $C^1$, tangent at $0$ to the eigenspace $E_0^c$ of the neutral eigenvalue $\lambda_k$.  A local center manifold is the graph of a $C^1$ function $\varphi_f: E_0^c \cap B' \rightarrow E_0^s\oplus E_0^u$, and is locally invariant, in the sense that $f(\Wloc)\cap B'\subset \Wloc$. The center manifold is not unique in general, but all center manifolds defined with respect to the ball $B'$ must contain the set of points which never escape from $B'$ under forward and backward iterations.  A weak uniqueness property can therefore be formulated as follows: if $f^{n}(x)\in B'$ for all $n\in\Z$, then $x\in \Wloc$. 

For the rest of the section, fix a local center manifold $\Wloc$ defined with respect to the ball $B'$. We show 
that the map $f$ restricted to $\Wloc$ is quasiconformally conjugate to an analytic map, in a two-step argument. Most of the analysis will be similar to Sections \ref{sec:structure} and \ref{sec:qc}, so we refer the reader to these sections for most proofs, and  we will only outline the differences, whenever they occur. We prove the following:

\begin{thmx}\label{thm:QC2} Let $f$ be a holomorphic germ of diffeomorphisms of $(\C^n,0)$. Suppose $df_0$ has eigenvalues $\lambda_j$, $1\leq j\leq n$, with $|\lambda_k|=1$ for some $k$ and $|\lambda_j|\neq 1$ when $j\neq k$.  Let $\WlocO$ be a $C^1$-smooth local center manifold of the fixed point $0$. There exist neighborhoods $W, W'$ of the origin inside $\WlocO$ such that $f:W\rightarrow W'$ is quasiconformally conjugate to a holomorphic diffeomorphism $h:(\Omega,0)\rightarrow (\Omega',0)$, $h(z) = \lambda_k z + \bigO(z^2)$, where $\Omega, \Omega'\subset\C$.

Moreover, the conjugacy map is holomorphic on the interior of $Z$ rel $\WlocO$, where $Z$ is the set of points that stay in $W$ under all forward and backward iterations of $f$.
\end{thmx}

\begin{remark}
Note that if $|\lambda_k|=1$ and $|\lambda_j|<1$ for all $j\neq k$ or $|\lambda_j|>1$ for all $j\neq k$, then the proof is identical to the proof of Theorem \ref{thm:QC}. 
\end{remark}

It is worth mentioning that the set $Z$ from Theorem \ref{thm:QC2} belongs to the intersection of all center manifolds defined relative to the ball $B'$. 
 
Denote by $J$ the standard almost complex structure of $\C^n$.
Consider a ball $B$ containing $0$, such that $\overline{B}\subset B'$.
We  first endow the two-dimensional real manifold $\Wloc$ with a $C^1$-smooth almost complex structure $J'$, induced by the restriction of the Riemannian metric of $\C^n$ to $\Wloc$. By integrating the almost complex structure $J'$, we show that the map $f$ on $W=\Wloc\cap B$ is conjugate to a map $g:U\subset \C \rightarrow \C$ of class $C^1$, via a  $(J',i)$-holomorphic  conjugacy map $\phi$, as in the diagram below:
\[
\diagup{W}{W'}{U}{U'}{f}{\phi}{\phi}{g}
\]

\smallskip
We will estimate how far  $g$ is from being an analytic map by measuring how far the tangent space $T_x\Wloc$ is from being a complex subspace of $T_x\C^n$, when $x\in W$. To carry on the analysis, 
we fix some notations for the dynamically relevant sets for $f$ and $g$. For each $n\geq 0$, let $W_n$ (respectively $W_{-n}$) be the set of points whose first $n$ backward (respectively forward) iterates remain in $W$. Similarly, we define the sets $U_n=\phi^{-1}(W_n)$ and $U_{-n}=\phi^{-1}(W_{-n})$ for the map $g$.

With these notations, $W_{\infty}$ and $W_{-\infty}$ represent the set of points from $W$ that do not escape from $B$ under backward, and respectively forward iterations. 
Let $Z := W_{\infty}\cap W_{-\infty}$. 
For simplicity, let $X=U_{\infty}$ and $Y=U_{-\infty}$.

\begin{prop}\label{horiz-2}$\ $
\begin{enumerate}
\item[a)] The tangent space $T_{x}\Wloc$ at any point $x\in Z$ is a complex line $E^c_x$ of $T_x\C^n$. The line field over $Z$ is $df$-invariant. 
\item[b)] There exists $\rho<1$ such that  for all integers $m,n\geq 0$ and for all $x\in W_n\cap W_{-m}$ and  $v\in T_{x}\Wloc$ the following estimate holds
\begin{equation*}
\ang \left(T_{x}\Wloc, {\rm Span}_{\C}\{v\}\right) = \bigO\left(\rho^{\min(m,n)}\right).
\end{equation*}
\end{enumerate}
\end{prop} 
\proof Part a) follows from the fact that for $x\in Z$
\[
T_x\Wloc= \left(\bigcap_{n\geq 0} df^n_{f^{-n}(x)} \mathcal{C}^{cu}_{f^{-n}(x)}\right) \cap \left(\bigcap_{n\geq 0} df^{-n}_{f^{n}(x)} \mathcal{C}^{cs}_{f^{n}(x)}\right),
\] 
and the counterpart of Proposition \ref{horiz}, which is straightforward. 
For part b) we observe that since $x\in W_n$, the tangent vector $v$ belongs to $df^n_{f^{-n}(x)}\mathcal{C}^{cu}_{f^{-n}(x)}$ which, by Remark \ref{rem:rho}, is a cone of angle opening $\alpha_1=\bigO(\rho^n)$ inside $\mathcal{C}^{cu}_{x}$, centered around $E^{cu}_x$. Similarly, since $x\in W_{-m}$, the tangent vector $v$ belongs to $df^{-m}_{f^{m}(x)}\mathcal{C}^{cs}_{f^{m}(x)}$, which is a cone of angle opening $\alpha_2=\bigO(\rho^m)$ inside $\mathcal{C}^{cs}_{x}$, centered around $E^{cs}_x$. Hence $v$ belongs to the complex cone centered around $E^c$, of angle less than the maximum of the angles $\alpha_1$ and $\alpha_2$. As in the proof of Proposition \ref{prop:span}, it follows that both $T_x\Wloc$ and $\mbox{Span}_{\C}\{v\}$ are included in this cone.
\qed

\begin{cor} \label{cor:submanifoldCn} Let $\mbox{int}^c(Z)$ denote the interior of $Z$ relative to $\Wloc$.  The set $\mbox{int}^c(Z)$ is a complex submanifold of $\C^n$ of complex dimension $1$. The conjugacy map $\phi : \mbox{int}(X\cap Y)\subset\C\rightarrow \mbox{int}^c(Z)\subset\C^n$ is holomorphic.
\end{cor}

\begin{lemma}\label{estimate:f2} There exists a constant $C$ such that for every $m,n\geq 1$, 
\[
\| \dbar_{J'} f\|_{W_n\cap W_{-m}} <C \rho^{\min(m,n)},
\]
where  $\dbar_{J'}f$ is the derivative of $f$ with respect to the almost complex structure $J'$ on $\Wloc$. 
\end{lemma}
 
The proof of this lemma uses Proposition \ref{horiz-2}. The argument is the same as in the proof of  Lemma \ref{estimate:f}, so we omit it here. 
 
\begin{prop}\label{prop:dbarCn} There exists a constant $C'$ such that for every $m,n\geq 1$, 
\[
| \dbar g(z)| <C' \rho^{\min(m,n)}, \ \ \mbox{for all}\ z\in U_n\cap U_{-m}.
\]
\end{prop}
\proof
The proof follows from Lemma \ref{estimate:f2} and Corollary \ref{cor:estimate-g2}.
\qed

In terms of Beltrami coefficients the proposition above implies that there exist $\rho<1$ and $M$ independent of $m, n$ such that 
\begin{equation}\label{eq:min-bound}
\sup_{z\in U_{n}\cap U_{-m}}|\mu_{g}(z)|<M\rho^{\min(m,n)}, \ \ \mbox{for all}\ m,n\geq 0.
\end{equation}

\begin{cor}\label{cor:dbar-XY} $\dbar g=0$ on $X\cap Y$.
\end{cor}

Note that in Section \ref{sec:structure} we obtained that the $\dbar$-derivative of $g$ is $0$ on the entire set $X$, whereas when we have stable, neutral and unstable eigenvalues we can only show that the  $\dbar$-derivative of $g$ is $0$ on $X\cap Y$.

Let $\sigma_0$ denote the standard almost complex structure of the plane, represented by the zero Beltrami differential. 
For $n\geq 0$, consider the Beltrami differential $\sigma_n$ on $U_n$, given by $\sigma_n = (g^{-n})^*\sigma_0$. Similarly, we define the Beltrami differentials $\sigma_{-n}$ on $U_{-n}$ by $\sigma_{-n} = (g^{n})^*\sigma_0$.

\begin{lemma}\label{lema:norm-inf}
There exists a constant $\kappa<1$, such that for all integers $n\geq0$ 
\[
\|\sigma_{n}\|_{\infty} = \|\sigma_{-n}\|_{\infty}<\kappa.
\]
\end{lemma}
\proof
Let  $z\in U_{-n}$. Let $z_{j} = g^j(z)$ for $0\leq j\leq n$ denote the $j$-th iterate of $z$ under the map $g$.
By Lemma \ref{prop:inclusions} part c), we know that  $z_{j}\in U_{j}\cap U_{-(n-j)}$ for all $0\leq j\leq n$. 
Hence, by Equation \eqref{eq:min-bound},
\[
\sup\limits_{z\in U_j \cap U_{-(n-j)}} |\mu_g(z)|<M\rho^{\min(j, n-j)}, \ \ 1\leq j\leq n.
\]

As in the proof of Lemma \ref{lema:bounded}, we show that the dilatation $K(g^n, z)$ is bounded by a constant independent of $n$ and the choice of $z$. We have 
\begin{equation}\label{eq:Kfn2}
K(g^n, z)\leq \prod_{j=0}^{n-1} K(g, z_{j}),
\end{equation}
where $K(g,z_{j}) =1+\bigO\left(\rho^{\min(j, n-j)}\right)$ and the conclusion follows.
\qed

We now use the estimates obtained on $\dbar g$ to prove the following theorem. 

\begin{thm}\label{thm:conformal2} 
The map $g^{-1}:U_1\rightarrow U_{-1}$ is quasiconformally conjugate to an analytic map.
\end{thm}
\proof Consider the measurable function $\mu:U \rightarrow \C$, given by 
\begin{equation*}\label{eq:mu2}
\mu = \left\{\begin{array}{ll}
                 \sigma_{n}& \mbox{on}\  U_n-U_{n+1},\ \mbox{for}\ n\geq 0 \\
		\sigma_{-n}& \mbox{on}\ X\cap \left(U_{-n}-U_{-(n+1)} \right),\ \mbox{for}\ n\geq 0 \\
		\sigma_{0}& \mbox{on}\  X\cap Y.
	\end{array}\right.
\end{equation*}
Then $\|\mu\|_{\infty}<1$ by Lemma \ref{lema:norm-inf}. Thus $\mu$ is a Beltrami coefficient, which is $g^{-1}$ invariant by construction, {\it i.e.} $(g^{-1})^{*}\mu = \mu$ on $U_1$. The set $X\cap Y$ is forward and backward invariant so $(g^{-1})^{*}\sigma_0 = \sigma_0$ on $X\cap Y$, by Corollary \ref{cor:dbar-XY}. The invariance of $\mu$ on $U-X$ is discussed in the proof of Theorem \ref{thm:conformal}. The only new case to check is when $n>0$ and  $z\in X\cap \left(U_{-n}-U_{-(n+1)} \right)$. By Lemma \ref{prop:inclusions},
\[
g\left(X\cap \left(U_{-n}-U_{-(n+1)} \right)\right)\subset X\cap \left(U_{-(n-1)}-U_{-n} \right).
\]
We have the following sequence of equalities
\[
\left(g^{-1}\right)^{*}\sigma_{-n}(z)=\left(g^{-1}\right)^{*}\left(g^{n}\right)^*\sigma_0(z)=\sigma_{-(n-1)}(z),
\]
which shows that $\left(g^{-1}\right)^{*}\mu = \mu$ on $X-Y$ as well.

The Measurable Riemann Mapping Theorem concludes the proof.
\qed

We proved that $f:W\rightarrow W'$ is quasiconformally conjugate to a holomorphic diffeomorphism $h:(\Omega,0)\rightarrow (\Omega',0)$,  where $\Omega, \Omega'\subset\C$. By Corollary \ref{cor:submanifoldCn} and Theorem \ref{thm:conformal2} it follows that the quasiconformal conjugacy map is holomorphic on the interior of $Z$, the set of points that remain in $W$ under all forward and backward iterates of $f$. The fact that $h(z) = \lambda_k z + \bigO(z^2)$, where $\lambda_k$ is the neutral eigenvalue of $df_0$, follows from the generalization of Naishul's theorem due to Gambaudo, Le Calvez, and P\'ecou \cite{GLP}. This concludes the proof of Theorem \ref{thm:QC2}.

\medskip

\noindent{\bf Types of hedgehogs.} Let $\mathcal{H}$ denote the connected component
containing $0$ of the set $Z$. Then $\mathcal{H}$ is the maximal hedgehog associated to the neighborhood $B$ of the origin. Using Theorem \ref{thm:QC2} and the local dynamics of the holomorphic germ $h$ of $(\C,0)$ with an indifferent fixed point at $0$, we can further describe the dynamical nature of the hedgehog. 

If $\lambda_k$ is a root of unity, $\lambda_k=e^{2\pi i p/q}$, and the parabolic multiplicity of $h$ at $0$ is $\nu$, then the hedgehog $\mathcal{H}$ of $f$ is the closure of $2\nu q$ holomorphic petals $\mathcal{P}_{inv}$, which are invariant under $f^{q}$ and $f^{-q}$, and where points converge both forward and backward to $0$. In addition, when $k=n$ in Equation \eqref{eq:k}, we can prove as in Theorem \ref{thm:F} the existence of holomorphic one-dimensional repelling petals with holomorphic outgoing Fatou coordinates. 
When $k\neq 1,n$, we can fix any center manifold $\Wloc\supset \mathcal{H}$ and use the quasiconformal conjugacy to construct $\nu q$ one-dimensional attracting and repelling petals in $\Wloc$, with the same regularity as the center manifold, consisting of points whose forward, respectively backward, orbit is contained in $\Wloc$ and converges to 0. However, the attracting/repelling petals will change as we change the center manifold. To visualize the phenomenon better, one may think of slicing the parabolic-attracting basin of $0$ (of complex dimension $k$), and the parabolic-repelling basin of $0$ (of dimension $n-k+1$) with different center manifolds. Theorem \ref{thm:QC2} also guarantees the existence of  quasiconformal incoming/outgoing Fatou coordinates $\varphi^i/\varphi^o$, with holomorphic transition maps $\varphi^i\circ (\varphi^{o})^{-1}:\mathcal{P}_{inv}\rightarrow\C$. 

If $\lambda_k=e^{2\pi i \alpha}$,  $\alpha\notin \Q$, and $h$ is linearizable at $0$ (that is, analytically conjugate to the rigid rotation $z\rightarrow \lambda_k z$ in a neighborhood of $0$ in $\C$), then $\mathcal{H}$ contains a holomorphic disk, called a Siegel hedgehog in our context. Lastly, if the angle $\alpha$ is irrational and $h$ is not linearizable at the origin, then $\mathcal{H}$ is a Cremer hedgehog, with a complicated topology: $\mathcal{H}$ has no interior, and is non-locally connected at any point different from the origin.

\section{Semi-parabolic germs of $(\C^n,0)$}\label{sec:para-Cn}

Let $f$ be a germ of diffeomorphisms of $(\C^n,0)$, with an isolated fixed point at the origin, with a neutral eigenvalue $\lambda=e^{2 \pi i p/q}$ a root of unity, $k-1$ eigenvalues $\lambda_2,\ldots, \lambda_k$ inside the unit disk, and $n-k$ eigenvalues $\lambda_{k+1},\ldots, \lambda_n$ outside of the unit disk. Let $\nu q +1$ be the order of $0$ as a fixed point of the equation $f^{q}(x)=x$. The germ $f(x,y,z) = (x_1,y_1,z_1)$ may be written near the origin in the form 
\begin{eqnarray*}
x_1 &=& \lambda x + f_1(x,y,z) \\
y_1 &=& A_2 y + f_2(x,y,z) \\
z_1 &=& A_3 z+f_3(x,y,z)
\end{eqnarray*}
where $x\in\C$, $y\in\C^{k-1}$, $z\in \C^{n-k}$, $A_2$ and $A_3$ are Jordan blocks corresponding to the attracting/repelling eigenvalues, and $f_i$, $1\leq i\leq 3$ are holomorphic functions which vanish at the origin along with their first order derivatives.

Let $r> 2\nu q+1$ be a positive integer. Consider as before a ball $B'$ containing the origin on which $f$ is partially hyperbolic. By choosing a sufficiently small neighborhood $B$ of the origin, compactly contained in $B'$, we may apply the Center Manifold Theorem to obtain unique analytic local strong stable/strong unstable manifolds $W_{\rm loc}^{ss}$, $W_{\rm loc}^{uu}$, and (non unique) center-stable, center-unstable and center manifolds, $W_{\rm loc}^{cs}$, $W_{\rm loc}^{cu}$, $W_{\rm loc}^{c}$ of class $C^r$, embedded in $B'$ and tangent to the linear subspaces $E^s_0$, $E^u_0$, $E^{cs}_0$, $E^{cu}_0$ and $E^c_{0}$. The center-stable manifold $W_{\rm loc}^{cs}=\{z=\varphi(x,y)\}$ has real dimension $2k$ and is the graph of a $C^r$-smooth function $\varphi:E^{cs}_0\cap B\rightarrow E^u_0$ with $\varphi(0,0)=0$. The center-unstable manifold $W_{\rm loc}^{cu}=\{y=\psi(x,z)\}$ has real dimension $2(n-k+1)$ and is the graph of a $C^r$-smooth function $\psi:E^{cu}_0\cap B\rightarrow E^s_0$, with $\psi(0,0)=0$. The center-stable and center-unstable manifolds intersect transversely, and their intersection is the center manifold $W_{\rm loc}^c$.

Let $\Lambda^{-}$ (respectively $\Lambda^{+}$) be the set of points that do not escape from $B$ under forward (respectively backward) iterations. Clearly $\Lambda^{-}\subset W_{\rm loc}^{cs}$ and $\Lambda^{+}\subset W_{\rm loc}^{cu}$.

We define the local parabolic-attracting (respectively parabolic-repelling) basins of $0$ to be the set of points from $\Lambda^{+}$ (respectively $\Lambda^{-}$) which converge to $0$ under forward (respectively backward) iterations, locally uniformly in $W_{\rm loc}^{cs}$ (respectively in $W_{\rm loc}^{cu}$). 

We will construct holomorphic Fatou coordinates on the parabolic basins. In the case when we have one neutral eigenvalue and all the other eigenvalues are inside the unit disk ($k=n$), this is done in \cite{U1} and \cite{Ha}. When we allow eigenvalues both inside and outside the unit disk we cannot apply the results in \cite{U1}, \cite{Ha} directly to the germ $f$, but we will be able to apply them to the germ $f$ restricted to the center-stable, and respectively to the center-unstable manifold of $0$. Even if $W_{\rm loc}^{cs}$ and $W_{\rm loc}^{cu}$ are only $C^r$-smooth, the parabolic basins are holomorphic, since they are contained in the analytic parts $\Lambda^{-}/\Lambda^{+}$ of the center-stable/center-unstable manifolds. 

\begin{prop}\label{horiz2} 
Let $W^{cs}_{\rm loc}$ and $W^{cu}_{\rm loc}$ be any local center-stable and center-unstable  manifolds of $0$ relative to $B'$. 
\begin{itemize}
\item[a)] The tangent space $T_{x}W^{cs}_{\rm loc}$ at any point $x\in \Lambda^{+}$ is a complex subspace $E^{cs}_x$ of $T_x\C^n$. The  distribution  over $\Lambda^{+}$ is $df^{-1}$-invariant, {\it i.e.} $df^{-1}_x(E^{cs}_x)=E^{cs}_{f^{-1}(x)}$ for every point $x\in \Lambda^{+}$.
\item[b)] The tangent space $T_{x}W^{cu}_{\rm loc}$ at any point $x\in \Lambda^{-}$ is a complex subspace $E^{cu}_x$ of $T_x\C^n$. The  distribution over $\Lambda^{-}$ is $df$-invariant, {\it i.e.} $df_x(E^{cu}_x)=E^{cu}_{f(x)}$ for every point $x\in \Lambda^{-}$.
\end{itemize}
\end{prop}
\proof We have:
\[
E^{cu}_x= \bigcap_{n\geq 0} df^n_{f^{-n}(x)} \mathcal{C}^{cu}_{f^{-n}(x)}\ \ \forall x\in \Lambda^{+} \ \mbox{ and }\  E^{cs}_x= \bigcap_{n\geq 0} df^{-n}_{f^{n}(x)} \mathcal{C}^{cs}_{f^{n}(x)}\ \ \forall x\in \Lambda^{-}.
\] 
The proof follows as in Proposition \ref{horiz}.
\qed

We outline here the construction of the Fatou coordinate for the parabolic-attracting basin, as the case of the parabolic-repelling basin is obtained analogously, by replacing $f$ by $f^{-1}$.

First we do a holomorphic change of coordinates in a neighborhood of the origin  to straighten the analytic strong stable manifold $W_{\rm loc}^{ss}$ of the fixed point $0$ to the axis $\{x=0, z=0\}$ and the analytic strong unstable manifold $W_{\rm loc}^{uu}$ to the axis $\{x=0, y=0\}$. 
Next we straighten the center-stable manifold $W^{cs}_{\rm loc}$ to the plane $z=0$, and the center-unstable manifold $W^{cu}_{\rm loc}$ to the plane $y=0$ via a local change of coordinates  $\phi(x,y,z)=(x_1,y_1, z_1)$, $x_1=x$, $y_1=y-\psi(x,z)$, and $z_1=z-\varphi(x,y)$,  of class $C^r$, which leaves $W_{\rm loc}^{ss}$ and $W_{\rm loc}^{uu}$ invariant. In the new coordinates, the center manifold $W_{\rm loc}^c$ is just the $x$-axis. The map $\phi$ is in fact analytic where $\psi$ and $\varphi$ are analytic. 

It is easy to see that $\phi(0,y,0)=(0,y,0)$ and  $\phi(0,0,z)=(0,0,z)$, 
hence in the new coordinate system we can think of $\pi_2:W_{\rm loc}^{ss}\rightarrow \C^{k-1}$, $\pi_2(0,y,0)=y$ as a natural complex coordinate on $W_{\rm loc}^{ss}$ and of $\pi_3:W_{\rm loc}^{ss}\rightarrow \C^{n-k}$, $\pi_3(0,0,z)=z$, as a natural complex coordinate on $W_{\rm loc}^{uu}$. 
The center manifold is only of class $C^r$, but we can put a reference complex structure on the center manifold to make $\pi_1:W_{\rm loc}^c\rightarrow \C$, $\pi_1(x,0,0)=x$ 
a complex coordinate on $W_{\rm loc}^c$. Consider a new complex structure on $\C^n$ by taking the product of the complex manifolds $\C\times\C^{k-1}\times \C^{n-k}$, each considered with the complex structures described above.  The new almost complex structure on $\C^n$ agrees with the usual almost complex structure of $\C^{n}$ up to order $r$.

Consider now the restriction of the map $f$ to the center-stable manifold. Let ${\rm int}^{cs}(\Lambda^-)$ denote the interior of the set $\Lambda^-$ rel $W_{\rm loc}^{cs}$.  Notice that  $f|_{W^{cs}_{\rm loc}}$ is holomorphic on  ${\rm int}^{cs}(\Lambda^-)$. In the new coordinates, $f|_{W^{cs}_{\rm loc}}$ has the form $f(x,y)=(x_{1},y_{1})$, with
\begin{equation}\label{eq:NF}
\left\{\begin{array}{l}
    x_{1}= \tilde{f}_1(x,y)\\
    y_{1}= g(y) + \tilde{f}_2(x,y),
    \end{array}\right.
\end{equation} 
where $\tilde{f}_1$, $\tilde{f}_2$ are germs of $C^r$-smooth functions of $(\C^k,0)$, and $g$ is a germ of holomorphic transformation of $(\C^{k-1},0)$. The restriction $f|_{W^{ss}_{\rm loc}}$ is analytically conjugate to $g$.  The map $g$ is contracting with $dg_0=A_2$, the Jordan block corresponding to the attracting eigenvalues $\{\lambda_j\}_{2\leq j\leq k}$. 
The germs $\tilde{f}_1$ and $\tilde{f}_2$ are holomorphic on ${\rm int}^{cs}(\Lambda^-)$. Moreover, $\tilde{f}_1(0,y)=0$, $\tilde{f}_2(0,y)=0$, 
and $\partial_x\tilde{f}_1$ at $(0,0)$ is $\lambda$.

We consider the Taylor series expansion of $\tilde{f}_1$ around a point $(0,y)\in W^{ss}_{\rm loc}$. 
The map  $\tilde{f}_1$ is not holomorphic, so we think of it as a $C^r$-smooth function in the variables $x$, $\bar{x}$ and $y$ and use the Wirtinger derivatives. Consider now an integer $1< r'<r$. We obtain:
\begin{equation}\label{eq:T1}
\tilde{f}_1(x, y) = \sum\limits_{1 \leq i + j \leq r'} \frac{\partial_x^i\partial_{\bar{x}}^j \tilde{f}_1(0, y) }{i!j!}x^i\bar{x}^j + \bigO_y(|x|^{r'+1}). 
\end{equation}
The tail $\bigO_y(|x|^{r'+1})$ is a $C^{r}$-smooth function of $x$ and $\bar{x}$. It is a $C^{r-r'}$-smooth function of $y$.

The map $\tilde{f}_1$ is holomorphic on ${\rm int}^{cs}(\Lambda^-)$. 
Since $W^{ss}_{\rm loc} = \{(0, y)\}$ lies in the closure of ${\rm int}^{cs}(\Lambda^-)$, it follows by continuity  that all  partial derivatives of $\bar{\partial}\tilde{f}_1$ vanish up to order $r'$ at $(0,y)$. Therefore the Taylor series expansion \eqref{eq:T1} reduces to 
\begin{equation}\label{eq:T2}
\tilde{f}_1(x, y) =  \sum\limits_{i=1}^{r'}  a_i(y) x^i + \bigO_y(|x|^{r'+1}),
\end{equation}
where the coefficients $ a_i(\cdot)$, $1\leq i\leq r'$, are $C^{r-i}$-smooth functions from $(\C^{k-1}, 0)$ to $\C$, with $a_1(0)=\lambda$.  Note that these coefficients are in fact holomorphic on ${\rm int}^{cs}(\Lambda^-)$. As in \cite{U1}, \cite{Ha}, \cite{RT}, we will do a series of changes of variables in the variable $x$ to make the coefficients $a_i(\cdot)$ constants, for $1\leq i\leq  r'$. The variable $y$ will remain unchanged.  We first make $a_1(y)=\lambda$ by considering the coordinate change $X = u(y)x$, $Y=y$, where $u$ is a germ of a $C^{r-1}$ function on $(\C^{k-1},0)$. We need to find $u$ so that 
\begin{eqnarray*}
X_{1}&=& u(y_{1})x_{1}=u\left(g(y) +\tilde{f}_2(x,y)\right)\left(  a_{1}(y)x + a_2(y)x^2 +\ldots\right)\\
    &=& u\left(g(Y) + \ldots \right) \left( a_{1}(Y)X/u(Y) +\ldots\right)=\frac{u(g(Y))a_{1}(Y)}{u(Y)}X + \ldots
\end{eqnarray*}
We seek a function $u$ such that $\frac{u(g(Y))a_{1}(Y)}{u(Y)}=\lambda$. We successively substitute $Y$ with $g(Y)$ in this equation and obtain 
\begin{equation}\label{eq:u(Y)}
u(Y) = \prod\limits_{n=0}^{\infty} \frac{a_{1}\left(g^{n}(Y)\right)}{\lambda},
\end{equation}
which converges in a neighborhood of $0$ since $\|g\|<1$ and $a_{1}(0) = \lambda$.

Suppose now that the first $i-1$ coefficients are constants and proceed by induction on $i$. The base case was treated above. 

Consider the coordinate transformation $X=x+v(y)x^{i}$, $Y=y$, where $v$ is a germ of a $C^{r-i}$ function on 
$(\C^{k-1},0)$. We have
\begin{eqnarray*}
X_{1}&=& x_{1} + v(y_{1})x_{1}^{i} = \sum\limits_{j=1}^{i-1} a_{j}x^{j} + \left(a_{i}(y)+\lambda^i v(g(y))\right)x^{i}+\ldots\\
    &=& \sum\limits_{j=1}^{i-1}a_{j}X^{j} +\left(a_{i}(Y) + \lambda^i v(g(Y)) - \lambda v(Y)\right)X^{i} + \ldots
\end{eqnarray*}

We choose $v$ such that the coefficient of $X^{i}$ is constant. This yields the functional equation $\lambda v(Y)-\lambda^i v(g(Y))=a_{i}(Y)-a_{i}(0)$, with solution  
\begin{equation}\label{eq:v(Y)}
 v(Y) =\frac{1}{\lambda} \sum\limits_{j=0}^{\infty} \left( a_{i}\left( g^{j}(Y)\right)-a_{i}(0) \right)\lambda ^{j(i-1)}.
\end{equation}
The series converges in a neighborhood of $0$ since $\|g\|<1$. We proved the following result.

\begin{prop}\label{prop:NF1}
There exists a $C^{r-r'}$ local change of coordinates in which $f|_{W^{cs}_{\rm loc}}$ has the form $f(x,y)=(x_{1},y_{1})$, where
\begin{equation}\label{eq:NF}
\left\{\begin{array}{l}
    x_{1}= \lambda x+ \sum\limits_{i=2}^{r'}a_{i}x^{i} + \bigO_y(|x|^{r'+1}) \\
    y_{1}=  g(y) + \tilde{f}_2(x,y)
    \end{array}\right.
\end{equation} 
and $a_i$ are constants, $2\leq i\leq r'$. The change of coordinates is analytic on  ${\rm int}^{cs}(\Lambda^-)$.
\end{prop}

\begin{prop}\label{prop:NF2}
There exists a $C^{r-2\nu q -1}$ local change of coordinates in which $f|_{W^{cs}_{\rm loc}}$ has the form $f(x,y)=(x_{1},y_{1})$, where
\begin{equation}\label{eq:NF2}
\left\{\begin{array}{l}
    x_{1}= \lambda x+x^{\nu q+1} + \bigO_y(|x|^{2\nu q +1}) \\
    y_{1}=  g(y) + \tilde{f}_2(x,y)
    \end{array}\right.
\end{equation} 
The change of coordinates is analytic on  ${\rm int}^{cs}(\Lambda^-)$.
\end{prop}
\proof We start from the normal form in Equation \eqref{eq:NF} and perform the same coordinate changes as in \cite{U1}, \cite{Ha}, \cite{RT}. 
If $k$ is not congruent to $1$ modulo $q$ (i.e. $\lambda^{k}\neq \lambda$), then we can set
\[b = \frac{a_{k}}{\lambda-\lambda^{k}}
\]
and consider the coordinate transformation  $X=x+bx^{k}, Y=y$ to eliminate the term $a_{k}x^{k}$. The first term which we cannot eliminate in this way is $a_{\nu q+1}x^{\nu q +1}$. We can further make the coefficient $a_{\nu q+1}=\lambda$ by considering a transformation of the form $X=Ax, Y=y$, where $A$ is a constant such that $A^{\nu q}=a_{\nu q+1}$. One can also eliminate all the intermediate terms of degree $jq+1$, for $\nu<j<2\nu$, by doing changes of variables of the form $X=x+bx^{(j-\nu)q+1}, Y=y$, where $b=\frac{a_{jq+1}}{(2\nu-j)q}$.

The changes of coordinates are holomorphic on ${\rm int}^{cs}(\Lambda^-)$, so in particular they are holomorphic on the parabolic-attracting basin of $0$.
\qed

\noindent{\bf Holomorphic Fatou coordinates.}
The choice of $r$ in the beginning of this section guarantees that the local change of coordinates in Proposition \ref{prop:NF2} is at least $C^{1}$ smooth. 
We can use the $C^1$-normal form \eqref{eq:NF2} to show the existence of $\nu q$ disjoint domains attracted to $0$ and invariant under $f^{ q}$, as in \cite{Ha}, \cite{U1}. We call these domains parabolic-attracting petals and denote them by $\paratt^{j}(0)$, $1\leq j\leq \nu q$. The key observation is that the parabolic-attracting petals belong to the analytic part of the center-stable manifold and the restriction of the map $f$ to any $\paratt^{j}(0)$ is holomorphic. Therefore, by Theorem \ref{prop:NF2}, $f$ restricted to $\paratt^{j}(0)$ is analytically conjugate to the normal form \eqref{eq:NF2}. Assuming \eqref{eq:NF2},
we can then use the same analytic changes of variables as in \cite{Ha}, \cite{U1} on each $\paratt^{j}(0)$ to show the existence of a holomorphic incoming Fatou coordinate $\varphi^i:\paratt^{j}(0)\rightarrow \C$, verifying the functional equation $\varphi^i(f^{q}(p))=\varphi^i(p)+1$. The parabolic-attracting basin of $0$ has $\nu q$ components $\Bparatt^{j}(0)$,  $1\leq j\leq \nu q$,  which are defined as 
\[
\Bparatt^{j}(0) = \bigcup_{n=1}^{\infty}f^{-nq}\left(\paratt^{j}(0)\right).
\]
By putting $\varphi^i(p)=\varphi^i(f^{n}(p))-n$ on $f^{-nq}\left(\paratt^{j}(0)\right)$ we can extend  $\varphi^i$ to the entire set $\Bparatt^{j}(0)$. Analogously, we can construct holomorphic outgoing Fatou coordinates on each component of the parabolic-repelling basin of $0$.

Having holomorphic Fatou coordinates has important consequence for the global dynamics. Hakim \cite{Ha} and Ueda \cite{U1} have shown that in the case of polynomial automorphisms of $\C^2$ with a semi-parabolic fixed point at the origin, each component of the parabolic basin of $0$ is a Fatou-Bieberbach domain, {\it i.e.} a proper subset of $\C^2$ biholomorphic to $\C^2$. The main strategy is to show that each component of the parabolic basin is a trivial fiber bundle over $\C$ with fibers biholomorphic to $\C$, given by the level sets of the incoming Fatou coordinate.   
The same technique, together with the holomorphic incoming/outgoing Fatou coordinates discussed in this section yield the following more general result:

\begin{thmx}\label{thm:Cn} Let $n\geq 3$ and consider $f$  a dissipative polynomial diffeomorphism of $\C^n$ with a fixed point at $0$. Suppose $df_0$ has a neutral eigenvalue $\lambda=e^{2\pi i p/q}$, a stable eigenvalue $0<|\lambda_{2}|<1$, while the rest of the eigenvalues $\lambda_{3},\ldots, \lambda_{n}$ are outside of the closed unit disk.  
Suppose the equation $f^{q}(x)=x$ has multiplicity $\nu q+1$ at the origin. 
Then the parabolic-attracting basin of $0$ has $\nu q$ components and each component is biholomorphic to $\C^{2}$. 
\end{thmx}

\begin{corx}\label{cor:Cn} If $n=3$, then the parabolic-repelling basin of $0$ has $\nu q$ components and each component is biholomorphic to $\C^{2}$. 
\end{corx}


\end{document}